\documentclass[a4paper,12pt]{amsart}
\usepackage{cite}
\usepackage{amsthm}
\usepackage{amssymb}
\usepackage{amsfonts}
\usepackage{amsmath}
\usepackage{url}
\usepackage{graphicx}
\usepackage{tikz}
\def\R{\mathbb R}
\newtheorem{thm}{Theorem}[section]
\newtheorem{lemm}{Lemma}[section]
\numberwithin{equation}{section}
\newtheorem{cor}{Corollary}[section]

\begin{document}

\author{J. Nahas}
\address[J. Nahas]{\'Ecole Polytechnique F\'ed\'erale de Lausanne\\
MA B1 487\\ 
CH-1015 Lausanne}
\email{joules.nahas@epfl.ch}

\title{A decay property of solutions to the k-generalized KdV equation}
%\author{Joules Nahas }
\maketitle

\begin{abstract}
\noindent We use a Leibnitz rule type inequality for fractional
derivatives to prove conditions under which a solution $u(x,t)$ of the k-generalized KdV equation
is in the space $L^2(|x|^{2s}\,dx)$ for $s \in \mathbb R_{+}$.
 \end{abstract}

\section{Introduction}
\normalsize The the initial value problem for the modified
Korteweg-de Vries equation (mKdV),
\begin{align}
\partial_tu + \partial_x^3u+\partial_x(u^3)=0, \label{mkdv} \\
u(x,0)=u_0(x), \notag
\end{align}
has applications to fluid dynamics (see \cite{2009ChPhB..18.4074L}, \cite{1994JNS.....4..355R}), and plasmas (see \cite{PRUD}). It is also an
example of an integrable system (see \cite{PhysRevLett.19.1095}). 
Ginibre and Y. Tsutsumi  in \cite{g} proved
well-posedness in a weighted $L^2$ space.
In \cite{KPV1}, Kenig, Ponce, and
Vega proved local well-posedness for $u_0$ in the Sobolev space $H^s$, when
$s \ge \frac{1}{4}$ by a contraction mapping argument in mixed
$L_x^p$ and $L_T^q$ spaces.
Christ, Colliander, and Tao in \cite{MR2018661}
showed that \eqref{mkdv} was locally well-posed for $u_0 \in H^s$, when $s
\ge \frac{1}{4}$, by using a contraction mapping argument in the
Bourgain spaces $X_{s,b}$.
  Colliander, Keel,
Staffilani, Takaoka, and Tao proved global well-posedness for real
initial data $u_0 \in H^{s}$, $s > \frac{1}{4}$ in \cite{CKSTT}. Kishimoto in \cite{Kish} and Guo in \cite{MR2531556} proved global well-posedness for real data in the case $s=\frac{1}{4}$.

The focus of this work will be \eqref{mkdv}, but we will also consider the generalized Korteweg-de Vries equation,
\begin{equation}
\left\{
\begin{array}{c l}
  & \partial_tu + \partial_x^3u + \partial_x (u^{k+1})=0, \label{gkdv} \\
& u(x,0)=u_0(x),\textrm{ } x \in \mathbb R.
\end{array}
\right.
\end{equation}
When $k \ge 4$, local well posedness was obtained for initial data $u_0 \in H^s$ with $s \ge \frac{k-4}{2k}$ in \cite{KPV1}
using a contraction mapping argument in mixed
$L_x^p$ and $L_T^q$ spaces. When $k=3$, the optimal local well posedness result was proven 
by Tao in \cite{MR2286393} for $u_0 \in H^s$ with $s \ge -\frac{1}{6}$
 by using Bourgain spaces $X_{s,b}$.

Kato in \cite{Ka} with energy estimates, and the fact that the
operator 
\begin{equation}
\Gamma_K \equiv x+3t\partial_x^2
\notag
\end{equation}
commutes with
$\partial_t+\partial_x^3$, was able to prove the following: if $u_0
\in H^{2k}$ and $|x|^ku_0 \in L^2$ where $k \in \mathbb Z^{+}$, then
for any other time $t$ when the solution exists, $|x|^ku(t) \in
L_x^2$. Using slightly different techniques, we will prove the following
theorem that extends this result slightly to $k \in \mathbb R_+$.
\begin{thm}
\label{weak-decay} 
Suppose the initial data $u_0$ satisfies
$|x|^su_0 \in L^2$, and $u_0 \in H^{2s+\varepsilon}$, for
$\varepsilon >0$. Then for any other time $t$, the solution $u(x,t)$
to \eqref{gkdv} satisfies $|x|^su(x,t) \in L^2$.

When $s \ge \frac{1}{2}$, the result holds for $\varepsilon=0$. Namely,
if $|x|^{s}u_0 \in L^2$, and $u_0 \in H^{2s}$, then for any
other time $t$, the solution $u(x,t)$ to \eqref{gkdv} satisfies
$|x|^{s}u(x,t) \in L^2$.
\end{thm}
 
Analogous results for the NLS were first
proved by Hayashi, Nakamitsu, and M. Tsutsumi in \cite{MR847012},
\cite{MR880978}, and \cite{MR987792}. They used the vector field
\begin{equation}
\Gamma_S = x+2it\nabla,
\label{nls-gamma}
\end{equation}
which commutes with the operator $\partial_t -i\Delta$, and a
contraction mapping argument to show that if $u_0 \in
L^2(|x|^{2m}\,dx) \cap H^m$, where $m \in \mathbb N$, then the
solution $u(x,t)$ at any other time is also in the space
$L^2(|x|^{2m}\,dx) \cap H^m$. These results were extended to the case when $m \in \mathbb R_+$
by the author and G. Ponce in \cite{NP}. The corresponding results for the Benjamin-Ono equation were obtained in
\cite{PonceFons} 
by G. Ponce and G. Fonseca.

Inspired by these persistence results we prove the following as our main result.
\begin{thm}
\label{main} If $u(x,t)$ is a solution of
\begin{equation}
\notag
\left\{
\begin{array}{c l}
  & \partial_tu + \partial_x^3u + \partial_x (u^{k+1})=0, \\
& u(x,0)=u_0(x),\textrm{ } x \in \mathbb R,
\end{array}
\right.
\end{equation}
such that $u_0 \in H^{s'} \cap L^2(|x|^{s}\,dx)$, where 
$s \in (0,s']$.
If $k=2$, and
$s'\ge
\frac{1}{4}$, then $u(\cdot,t) \in H^{s'} \cap
L^2(|x|^{s}\,dx)$ for all $t$ in the lifespan of $u$.

If $k \ge 4$, and
$s \ge \frac{k-4}{2k}$, then $u(\cdot,t) \in H^{s'} \cap
L^2(|x|^{s}\,dx)$ for all $t$ in the lifespan of $u$.
\end{thm}

We only prove this property the most interesting case, \eqref{mkdv}.
Note that the cases in \eqref{gkdv} when $k=1$ or $4$ are excluded from
Theorem \ref{main}. We require
our technique to be adapted to Bourgain spaces for these nonlinearities, which is an
interesting open question.

The difficulty in the case of fractional decay lies in the lack of
an operator $\Gamma$ that sufficiently describes the relation
between initial decay, and properties of the solution at another
time (such as \eqref{nls-gamma}). In order to solve this problem, we develop a Leibnitz rule
type inequality for fractional derivatives.

We need some notation to illustrate this idea. If $f$ is a complex
valued function on $\mathbb R$, we let $f^{\wedge}$ (or $\hat{f}$)
denote the Fourier transform of $f$, and $f^{\vee}$ the inverse
Fourier transform. For $\alpha \in \mathbb R$, the operator
$D_x^{\alpha}$ is defined as $(D_x^{\alpha}f(x))^{\wedge}(\xi)\equiv
|\xi|^{\alpha}f^{\wedge}(\xi)$. Let $U(t)f$ denote the solution
$u(x,t)$ to the linear part of $\eqref{mkdv}$, with $u(x,0)=f(x)$.
Choose $\eta \in C_0^{\infty}(\mathbb R)$ with
$\textrm{supp}(\eta)\subset [\frac{1}{2},2]$ so that
\begin{equation}
\sum_{N \in \mathbb Z} (\eta(\frac{x}{2^N})+\eta(-\frac{x}{2^N}))=1
\textrm{ for }x \ne 0.
\notag
\end{equation}
Define the operator $Q_N$ on a function $f$ as
\begin{equation}
Q_N(f) \equiv ((\eta(\frac{\xi}{2^N})+\eta(-\frac{\xi}{2^N}))\hat{f}(\xi))^{\vee}. \notag
\end{equation}
If $\|\cdot\|_Y$ is a norm on some space of functions, we recall
that
\begin{equation}
\|Q_N(f)\|_{Y l_N^p} \equiv \|(\sum_{N \in \mathbb Z}|Q_N(f)|^p)^{\frac{1}{p}}\|_Y\notag.
\end{equation}

Using Duhammel's principle, we can formulate the problem
\eqref{mkdv} as an integral equation.
\begin{equation}
u(x,t) = U(t)u_0 -\int_{0}^{t}U(t-t')\partial_x(u^3(x,t'))\,dt'.
\notag
\end{equation}
Using a Fourier transform, we can see how to commute an $x$ past
$U(t)$,
\begin{align}
xU(t)f
& =
(-i\partial_{\xi}(e^{it\xi^3}\hat{f}))^{\vee}
\notag \\
& =
(3t\xi^2e^{it\xi^3}\hat{f}-ie^{it\xi^3}\partial_{\xi}\hat{f})^{\vee}
\notag \\
& =
U(t)(3t\partial_x^2f+xf).
\notag
\end{align}
We would like to use a similar argument with $|x|^{\frac{1}{8}}$
replacing $x$, but this would require that $D_{\xi}^{\frac{1}{8}}$
obey a product rule. We develop in inequality in Lemma
\ref{my_prod_rule} that is similar enough to the product rule that
will allow this argument to work.

With Lemma \ref{my_prod_rule}, we will require that
\begin{equation}
\left \|D_{\xi}^{\frac{1}{8}}Q_N(\frac{e^{it\xi^3}}{(1+\xi^2)^{\frac{1}{8}}}) \right \|_{L_{\xi}^{\infty}l_N^1}
< \infty. \label{main_part}
\end{equation}

With less sophisticated techniques, we prove Theorem {weak-decay} in Section 2. We show \eqref{main_part} in Section 3, then prove our main result
in Section 4. The proof of Lemma \ref{my_prod_rule} is almost
identical to the proof of a classical Leibnitz rule inequality.
Because this proof requires different techniques than the rest of
the paper, we present it in Appendix A.

We use the following notation throughout the paper. We let $A \lesssim B$ mean that the quantity $A$ is less than or
equal to a fixed constant times the quantity $B$. Let 
$\langle x \rangle \equiv (1+x^2)^{\frac{1}{2}}$, and similarly,
$\langle D_x \rangle$.

\section{Weak Persistence Result}

Using some standard estimates, we prove Theorem \ref{weak-decay} which is a weaker persistence property for IVP for the gKdV equation for low regularity solutions, but holds for more values of $k$ in \eqref{gkdv} than
our main result.

Following an argument by Kato, we multiply \eqref{gkdv}
by $\phi(x)u(x,t)$ for some function $\phi(x)$, and integrating over
$x$ and $t$, we use integration by parts to obtain
\begin{align}
& \int_{\mathbb R}\phi(x)u^2(x,T)\,dx-\int_{\mathbb R}\phi(x)u^2(x,0)\,dx
-3\int_{[0,T]}\int_{\mathbb R}\phi'(x)(\partial_x u)^2\,dx\,dt \notag \\
& \quad
+\int_{[0,T]}\int_{\mathbb R}\phi'''(x)u^2\,dx\,dt
+ \frac{k+1}{k+2}\int_{[0,T]}\int_{\mathbb R}\phi'(x)u^{k+2}\,dx\,dt=0.
\label{weight_eq}
\end{align}
Equation \eqref{weight_eq}, along with the following two interpolation lemmas are the primary tools for the
weak persistence result, Theorem \ref{weak-decay}.

\begin{lemm}
\label{inter-reg-dec}
Let $a,b >0$, and $w(x)> \varepsilon >0$ a locally bounded function. Assume that $\langle D_x \rangle^a f \in L^2(\mathbb R)$
and $ w^b(x) f \in L^2(\mathbb R)$. Then for any $\theta \in (0,1)$
\begin{equation}
\|\langle D_x \rangle^{\theta a}( w^{(1 - \theta)b}(x) f)\|_2 
\lesssim 
\| w^{b}(x) f\|_2^{1-\theta}\|\langle D_x \rangle^af\|_2^{\theta}.
\notag
\end{equation}
\end{lemm}
\begin{proof}
This is an easy consequence of the Three Lines Lemma, and the fact that
 \begin{equation}
\|\langle D_x \rangle^{z a}( w^{(1 - z)b}(x)f)\|_2
\notag
\end{equation}
is an analytic function in $z$ for $\Re z \in (0,1)$, for a dense set of functions 
in the space $H^a \cap L^2( w^{2b}(x)\,dx)$.
\end{proof}

\begin{lemm}
\label{23a}
For a solution $u=u(x,t)$ of \eqref{gkdv},
\begin{equation}
 \|\partial_xu\|_{L^{\frac{1}{s+\frac{1}{2}\varepsilon}}_xL^2_T}\leq c_T \|u_0\|_{H^{2s+\varepsilon}}.
\end{equation}
\end{lemm}
\begin{proof}
Consider the function
$$
F(z)=\int_{-\infty}^{\infty}\int_0^T D_x^{r(z)}(U(t)u_0)\, \psi(x,z)\,f(t)\,dt dx,
$$
where
$$
r(z)=(1-z)(1+2s+\varepsilon)+z (2s+\varepsilon),\;\;\;\frac{1}{q(z)}=\frac{z}{2}+(1-z),\;\;\;q=\frac{2}{2-2s-\varepsilon},
$$
$$
\psi(x,z)=|g(x)|^{q/q(z)}\,\frac{g(x)}{|g(x)|},\;\;\;\;\text{with}\;\;\;\;\|g\|_{L^{1/(1-s-\frac{1}{2}\varepsilon)}_x}=\|f\|_{L^2([0,T])}=1,
$$
which is analytic for $\Re z \in (0,1)$. Using that
$$
\|\psi(\cdot,0+iy)\|_2=\|\psi(\cdot,1+iy)\|_1=1,
$$
one gets from $H^{2s+\varepsilon}$ persistence and the Kato smoothing effect that
\begin{equation}
\label{23}
\begin{aligned}
&  \|\partial_xU(t)u_0\|_{L^{\frac{1}{s+\frac{1}{2}\varepsilon}}_xL^2_T} \leq c \|D_xU(t)u_0\|_{L^{\frac{1}{s+\frac{1}{2}\varepsilon}}_xL^2_T}\\
& \quad
\leq c \,\sup_{y\in\R} \| D_x^{1+2s+\varepsilon+iy} U(t)u_0\|_{L^{\infty}_xL^2_T}^{1-2s-\varepsilon}\,\sup_{y\in\R}\|D^{2s+\varepsilon+iy}_xU(t)u_0\|_{L^2_xL^2_T}^{2s+\varepsilon}
\notag \\
& \quad
\leq c_T \| D_x^{2s+\varepsilon}U(t)u_0\|_{2}.
\notag
\end{aligned}
\end{equation}

Inserting the estimate \eqref{23} in the proof of the local well posedness for \eqref{gkdv}, the result follows.
\end{proof}

\begin{proof}[Proof of Theorem \ref{main}]
Let $\phi_N$ be a smooth function such that
\begin{equation}
\phi_N(x) = \left\{
\begin{array}{c l}
  \langle x \rangle^{2s} & \textrm{if $|x| \leq N$,} \\
(2N)^{2s} & \textrm{if $|x|>3N$}.
\end{array}
\right.\notag \notag
\end{equation}
Then from \eqref{weight_eq},
\begin{align}
& \int_{\mathbb R}\phi_N(x)u^2(x,T)\,dx-\int_{\mathbb R}\phi_N(x)u^2(x,0)\,dx
 =
\notag \\
& \quad 
3\int_{[0,T]}\int_{\mathbb R}\phi'_N(x)(\partial_x u)^2\,dx\,dt
-\int_{[0,T]}\int_{\mathbb R}\phi_N'''(x)u^2\,dx\,dt
\notag \\
& \quad - \frac{k+1}{k+2}\int_{[0,T]}\int_{\mathbb R}\phi'_N(x)u^{k+2}\,dx\,dt.
\label{weight_eq_2}
\end{align}
We only prove the result in the case where $s < 2$ of
the KdV equation, when $k=1$. Our main result, Theorem \ref{main},
is stronger when $k=2$, and $k \ge 4$, and the proof for $s \ge
2$ or $k=3$ is similar. 
We will use results from
\cite{MR1230283}, which state that the smoothing effects and
Strichartz estimates that hold for the linearized KdV and mKdV also
hold for the KdV.

The $\phi_N'''(x)u^{2}$ term in the right hand side of
\eqref{weight_eq_2} can be bounded by the fact that $\phi_N'''(x)
\lesssim 1$ independently of $N$ for $s \leq \frac{1}{2}$, and $L^2$
persistence:
\begin{equation}
\left | \int_{[0,T]}\int_{\mathbb R}\phi_N'''(x)u^2\,dx\,dt \right |
\lesssim
T\|u\|_2^2.
\label{missed-term} 
\end{equation}

The bounds on the other terms on the right hand side of \eqref{weight_eq_2} 
depend on whether $s < \frac{1}{2}$ or $s \ge \frac{1}{2}$. We first
give the proof of the result in the case that $s < \frac{1}{2}$ 

Since $|\phi_N'(x)| \lesssim
\langle x \rangle^{2s-1}$ independently of $N$, we can bound the
first term on the right hand side of \eqref{weight_eq_2} by
\begin{equation}
\left | \int_{[0,T]}\int_{\mathbb R}\phi'_N(x)(\partial_xu)^2\,dx \,dt \right |
\lesssim
\|\langle x \rangle^{s-\frac{1}{2}}\partial_xu\|_{L_x^2L_T^2}^2.
\label{1_rhs_weight}
\end{equation}
Using \eqref{1_rhs_weight}, Lemma \ref{23a}, and the H\"{o}lder
inequality,
\begin{align}
&  \left |\int_{[0,T]}\int_{\mathbb R}\phi'_N(x)(\partial_x u)^2\,dx\,dt \right | \notag \\
& \quad \lesssim
\|\langle x \rangle^{s-\frac{1}{2}}\|_{\frac{2}{1-2(s+\frac{1}{2}\varepsilon)}}
\|D_xu(x,t)\|_{L_x^{\frac{1}{s+\frac{1}{2}\varepsilon}}L_T^2}
< \infty. \label{weight-eq-first-rhs-term}
\end{align}

For the $\phi_N'(x)u^{k+2}$ term in the right hand side
of \eqref{weight_eq_2}. We can bound this term with the H\"{o}lder
inequality,
\begin{align}
\left |
\int_{[0,T]}\int_{\mathbb R}\phi'_N(x)u^{3}\,dx\,dt
\right |
& \lesssim
 \| \langle x\rangle^{2s-1}|u|^{3} \|_{L_T^1L_x^1}
\notag \\
& \leq
\|  u \|_{L_T^1L_x^{\infty}}
\| \langle x\rangle^{s-\frac{1}{2}}u \|_{L_T^{\infty}L_x^2}^2
\notag \\
& \leq
T^{\frac{5}{6}} \|  u \|_{L_T^6L_x^{\infty}}
\|u \|_{L_T^{\infty}L_x^2}^2.
\label{weight-eq-sec-rhs-term}
\end{align}
Since $s-\frac{1}{2} < 0$, \eqref{weight-eq-sec-rhs-term} is finite by
the Strichartz estimates in \cite{MR1230283}, and $L^2$ persistence.

It follows from \eqref{weight_eq_2} that
\begin{align}
\left |\int_{\mathbb R}(\phi_N(x)u^2(x,T))\,dx \right |
& \leq
\int_{\mathbb R}|\phi_N(x)u^2(x,0)|\,dx
\notag \\
& \quad +
3\int_{[0,T]}\int_{\mathbb R}|\phi'_N(x)(\partial_x u)^2|\,dx\,dt
\notag \\
& \quad +
\frac{2}{3}\int_{[0,T]}\int_{\mathbb R}|\phi'_N(x)u^{3}|\,dx\,dt
\notag \\
& \quad +
\int_{[0,T]}\int_{\mathbb R}|\phi_N'''(x)u^2|\,dx\,dt.
\notag
\end{align}
By $|x|^{s}u_0 \in L^2$, \eqref{weight-eq-first-rhs-term}, \eqref{missed-term}, and
\eqref{weight-eq-sec-rhs-term}, the result follows.

We now consider the case that $s \in [\frac{1}{2},1)$. For the first term on the right hand side of
\eqref{weight_eq_2}, we use Lemma \ref{inter-reg-dec}, and $H^{2s}$
persistence to obtain
\begin{align}
\left | \int_{[0,T]}\int_{\mathbb R}\phi'_N(x)(\partial_xu)^2\,dx \,dt \right |
& \lesssim
\|\partial_x u \langle \phi_N'(x) \rangle^{\frac{1}{2}} \|_{L_T^2L_x^2}^2 
\notag \\
& \lesssim
\|\partial_x (u\langle \phi_N'(x) \rangle^{\frac{1}{2}}) \|_{L_T^2L_x^2}^2 
\notag \\
& \quad +
\|u(\langle \phi_N'(x) \rangle^{\frac{1}{2}})' \|_{L_T^2L_x^2}^2 
\notag \\
& \lesssim
\|\frac{\partial_x}{\langle D_x \rangle} \langle D_x \rangle(u\langle \phi_N'(x) \rangle^{\frac{1}{2}}) \|_{L_T^2L_x^2}^2 
\notag \\
& \quad +
\|u\langle x \rangle ^{s-\frac{3}{2}} \|_{L_T^2L_x^2}^2 
\notag \\
& \lesssim
\| \langle D_x \rangle(u\langle \phi_N'(x) \rangle^{\frac{1}{2}}) \|_{L_T^2L_x^2}^2 
\notag \\
& \quad+
\|u\langle x \rangle ^{s-\frac{3}{2}} \|_{L_T^2L_x^2}^2 
\notag \\
& \lesssim
\| \langle D_x \rangle^{2s} u \|_{L_T^2L_x^2}^{\frac{1}{s}} 
\| (\langle \phi_N'(x) \rangle^{\frac{s}{2s-1}})u \|_{L_T^2L_x^2}^{2-\frac{1}{s}} 
\notag \\
& \quad +
\|u\langle x \rangle ^{s-\frac{3}{2}} \|_{L_T^2L_x^2}^2 
\notag
\end{align}
Since $\langle \phi_N'(x) \rangle^{\frac{s}{2s-1}} \lesssim \phi_N^{\frac{1}{2}}(x)$, it follows that
\begin{align}
\left | \int_{[0,T]}\int_{\mathbb R}\phi'_N(x)(\partial_xu)^2\,dx \,dt \right |
& \lesssim
\| \langle D_x \rangle^{2s} u \|_{L_T^2L_x^2}^{\frac{1}{s}} 
\| \phi_N^{\frac{1}{2}}(x)u \|_{L_T^2L_x^2}^{2-\frac{1}{s}} 
\notag \\
& \quad +
\|u\langle x \rangle ^{s-\frac{3}{2}} \|_{L_T^2L_x^2}^2 
\label{case-half}
\end{align}

For the $\phi_N'(x)u^{k+2}$ term,
\begin{align}
\left |
\int_{[0,T]}\int_{\mathbb R}\phi'_N(x)u^{3}\,dx\,dt
\right |
& \lesssim
 \| \langle x\rangle^{2s-1}|u|^{3} \|_{L_T^1L_x^1}
\notag \\
& \leq
\|  u \|_{L_T^1L_x^{\infty}}
\| \langle x\rangle^{s-\frac{1}{2}}u \|_{L_T^{\infty}L_x^2}^2
\notag \\
& \leq
T^{\frac{5}{6}} \|  u \|_{L_T^6L_x^{\infty}}
\| \langle x\rangle^{s-\frac{1}{2}}u \|_{L_T^{\infty}L_x^2}^2
\label{big-s-3}
\end{align}
The term in \eqref{big-s-3} is finite from the first part of the proof since $s-\frac{1}{2} < \frac{1}{2}$.

From \eqref{weight_eq_2}, \eqref{missed-term}, \eqref{big-s-3}, \eqref{case-half}, 
the fact that $\phi_N(x) \lesssim \langle x \rangle^{2s}$
and our assumption on $u(x,0)$,
\begin{align}
& \|\phi_N^{\frac{1}{2}}(x)u^2(x,T)\|_{L_x^2}^2
\lesssim \|\langle x \rangle^s u^2(x,0)\|_{L_x^2}^2
+
\|u\langle x \rangle ^{s-\frac{3}{2}} \|_{L_T^2L_x^2}^2 
\notag \\
& \quad +\| \langle x\rangle^{s-\frac{1}{2}}u \|_{L_T^{\infty}L_x^2}^2
  +
\| \langle D_x \rangle^{2s} u \|_{L_T^2L_x^2}^{\frac{1}{s}} 
(\int_0^T\| \phi_N^{\frac{1}{2}}(x)u(x,t) \|_{L_x^2}^2\,dt )^{1-1/2s}
\notag \\ 
\label{stay-on-target}
& \quad +T\|u\|_{L_x^2}^2
+ T^{\frac{5}{6}} \|  u \|_{L_T^6L_x^{\infty}}
\|u \|_{L_T^{\infty}L_x^2}^2.
\end{align}
The applicaion of Bihari's inequality (see \cite{MR0079154}) to \eqref{stay-on-target}
yields a bound on $\|\phi_N^{\frac{1}{2}}(x)u(x,T)\|_{2}$ that is independent of $N$.
By taking $N$ to infinity, the result follows.

\end{proof}

\section{Estimating a Derivative}
We begin our computation of \eqref{main_part}. We will show that by
scaling out the fractional derivative, it will suffice to bound
\begin{equation}
\left |Q_N(\frac{e^{it\xi^3}}{(1+\xi^2)^{\frac{1}{8}}}) \right |. \notag
\end{equation}
Since the operator $Q_N$ is convolution with a function whose
Fourier transform is very localized, we require estimates on
\begin{equation}
\int_{\mathbb
R}\varphi_{\omega}(\xi-z)\frac{e^{itz^3}}{(1+z^2)^{\frac{1}{8}}}\,dz,
\label{project}
\end{equation}
where $\varphi_{\omega}$ is a function whose Fourier transform has
support near $\omega$.

 We will use a contour integral argument. Because of
this, we require estimates on the analytic continuation of
$\varphi_{\omega}$. These are contained in the following lemma.
\begin{lemm}
\label{anal_cont}
 Let $\xi \in \mathbb R$, $z=x+yi$ for $x,y \in \mathbb R$, $\varphi(\xi)$ be a function so that
$\hat{\varphi}(x)$ is a smooth function with support in
$[\frac{1}{2},2]$, and for $\omega \in \mathbb R \setminus \{0\}$,
let $\varphi_{\omega}(\xi)$ be the function with Fourier transform
$\hat{\varphi}(\frac{x}{\omega})$. Then $\varphi_{\omega}$ is an
entire function that obeys the following estimates.
\begin{equation}
|\varphi_{\omega}((\xi-z))| \lesssim
\left\{
\begin{array}{c l}
 \frac{|e^{2\omega y}-e^{\frac{1}{2}\omega y}|}{{\omega}^2y|\xi-z|^2} & \textrm{if $y \ne 0$ and $x \ne \xi$,}\\
 \\
 \frac{1}{|\omega|(\xi-x)^2} & \textrm{if $y = 0$ and $x \ne \xi$.}
\end{array}
\right
. \notag
\end{equation}
\end{lemm}
\begin{proof}
That $\varphi_{\omega}$ is entire follows from the Paley-Wiener
theorem. Let $y \ne 0$.
 Since
$\hat{\varphi}$ is a smooth function with support in
$[\frac{1}{2},2]$, we integrate by parts to obtain
\begin{align}
\varphi_{\omega}(\xi-z) & = \int_{\mathbb R}\hat{\varphi}(\frac{\zeta}{\omega})\frac{1}{i(\xi-z)}\frac{d}{d\zeta}e^{i\zeta (\xi-z)}\,d\zeta \notag \\
& = -\int_{[\frac{\omega}{2},2\omega]}\frac{1}{\omega}\hat{\varphi}^{'}(\frac{\zeta}{\omega})\frac{1}{i(\xi-z)}e^{i\zeta (\xi-z)}\,d\zeta \notag\\
& = \int_{[\frac{\omega}{2},2\omega]}\frac{1}{\omega}\hat{\varphi}^{'}(\frac{\zeta}{\omega})\frac{1}{(\xi-z)^2}\frac{d}{d\zeta}e^{i\zeta (\xi-z)}\,d\zeta \notag\\
& = -\int_{[\frac{\omega}{2},2\omega]}\frac{1}{{\omega}^2}\hat{\varphi}^{''}(\frac{\zeta}{\omega})\frac{1}{(\xi-z)^2}e^{i\zeta (\xi-z)}\,d\zeta. \label{phi_parts}
\end{align}
From \eqref{phi_parts} we conclude that
\begin{align}
|\varphi_{\omega}(\xi-z)| & \leq \left |\int_{[\frac{\omega}{2},2\omega]}\frac{1}{{\omega}^2}\hat{\varphi}^{''}(\frac{\zeta}{\omega})\frac{1}{(\xi-z)^2}e^{i\zeta (\xi-z)}\,d\zeta \right | \notag \\
& \leq \int_{[\frac{\omega}{2},2\omega]}\frac{1}{{\omega}^2}|\hat{\varphi}^{''}(\frac{\zeta}{\omega})|\frac{1}{|\xi-z|^2}e^{\zeta y}\,d\zeta \notag \\
& \leq c_{\varphi}\frac{|e^{2\omega y}-e^{\frac{1}{2}\omega y}|}{{\omega}^2y|\xi-z|^2}. \notag
\end{align}
The case $y=0$ follows from taking the limit as $y \rightarrow 0$ of
the first estimate.
\end{proof}
From Lemma \ref{anal_cont}, we can infer the following about the
analyticity of the integrand in \eqref{project}.
\begin{cor}
\label{anal_cont_cor} For $\xi \in \mathbb R$, the function
\begin{equation}
\varphi_{\omega}(\xi-z)\frac{e^{itz^3}}{(1+z^2)^{\frac{1}{8}}}
\end{equation}
is analytic on $\mathbb C \setminus \{z: |\Im{z}| \ge 1,
\Re{z}=0\}$.
\end{cor}
The estimate in Lemma \ref{anal_cont} has good $x$ dependence away
from $\xi$. To estimate \eqref{project} near $z=\xi$, we use an
analytic continuation of the integrand and the Cauchy integral
theorem, which we now describe.

The function $\varphi_{\omega}$ oscillates with frequency near
$\omega$. For a fixed $z_0 \in \mathbb R$, we think of the function
$\exp(itz^3)$ as oscillating with frequency $tz^2_0$ near the value
$z_0$. For $z=\xi$ where $t\xi^2 \ll \omega$, the function $\varphi_{\omega}$ oscillates much faster than $\exp(itz^3)$, so Lemma \ref{anal_cont}
shows that analytic continuation of
\begin{equation}
\varphi_{\omega}(\xi-z)\frac{e^{itz^3}}{(1+z^2)^{\frac{1}{8}}}
\notag
\end{equation}
changes this rapid oscillation into
decay, which yields good $\omega$ dependence for \eqref{project}. To formalize this, we
make the following definition. Given $t>0$, and $\omega>0$, we say
that $\xi \in \mathbb R$ is \textbf{near} if
\begin{equation}
|\xi| \leq \frac{1}{10}\sqrt{\frac{\omega}{t}}. \notag
\end{equation}
Where the oscillation of $\exp(itz^3)$ is much larger than $\omega$,
an analytic continuation of $\exp(itz^3)$ has a similar property. We
say that $\xi \in \mathbb R$ is \textbf{far} if
\begin{equation}
|\xi| > 10\sqrt{\frac{\omega}{t}}. \notag
\end{equation}
In the intermediate case where the oscillation of $\exp(itz^3)$ is
comparable to $\omega$, analytic continuation does not help. This is
where the worst behavior of the estimate occurs. We say that $\xi
\in \mathbb R$ is \textbf{intermediate} if
\begin{equation}
\frac{1}{10}\sqrt{\frac{\omega}{t}} < |\xi| \leq 10\sqrt{\frac{\omega}{t}}. \notag
\end{equation}
These heuristics are formalized in Lemma \ref{convol}, then used to
estimate \eqref{main_part} in Lemma \ref{finite}. We require an
elementary integral estimate for Lemma \ref{convol}.

One expects that since $\sin t \approx t$, then
\begin{align}
 \int_{[0,\pi]}\frac{e^{a\sin{s}}-e^{b\sin{s}}}{\sin{s}}\,ds
&  \approx
\int_{[0,\pi]}\frac{e^{a s}-e^{b s}}{s}\,ds
\notag \\
& =
\int_{[0,\pi]}\frac{e^{a s}-e^{b s}}{as}a\,ds \notag \\
& =
\int_{[0,\pi]}\frac{e^{t}-e^{\frac{b}{a} t}}{t}\,dt. \notag
%\lesssim (\pi\frac{a}{b}-1)+1+\frac{b}{\pi a}e^{-\frac{\pi a}{b}}.
\end{align}
This is what the next lemma proves.

\begin{lemm}
\label{e_sin} Let $a<b<0$. Then
\begin{equation}
\left | \int_{[0,\pi]}\frac{e^{a\sin{s}}-e^{b\sin{s}}}{\sin{s}}\,ds\right |
\lesssim (\pi\frac{a}{b}-1)+1+\frac{b}{\pi a}e^{-\frac{\pi a}{b}}.
\notag
\end{equation}
\end{lemm}
\begin{proof}
By making the change of variable $r=-(s-\frac{\pi}{2})$, we have for
an arbitrary function $f$,
\begin{equation}
\int_{[\frac{\pi}{2},\pi]}f(\sin{s})\,ds=\int_{[0,\frac{\pi}{2}]}f(\sin{r})\,dr. \notag
\end{equation}
Therefore,
\begin{align}
\left | \int_{[0,\pi]}\frac{e^{b\sin{s}}-e^{a\sin{s}}}{\sin{s}}\,ds\right |
& =
 2\int_{[0,\frac{\pi}{2}]}\frac{e^{b\sin{s}}-e^{a\sin{s}}}{\sin{s}}\,ds.
\label{e_sin_1}
\end{align}
Notice that for $s \in [0,\frac{\pi}{2}]$, $\frac{2s}{\pi} \leq
\sin{s} \leq 2s$. We use this to bound \eqref{e_sin_1}.
\begin{align}
\int_{[0,\frac{\pi}{2}]}\frac{e^{b\sin{s}}-e^{a\sin{s}}}{\sin{s}}\,ds & \lesssim
\int_{[0,\frac{\pi}{2}]}\frac{e^{\frac{2b}{\pi}s}-e^{2as}}{s}\,ds
\notag \\
& =
\int_{[0,\frac{\pi}{2}]}\frac{e^{\frac{2b}{\pi}s}-e^{2as}}{\frac{2b}{\pi}s}\frac{2b}{\pi}\,ds
\notag \\
& =
\int_{[b,0]}\frac{e^{\frac{\pi a}{b}r}-e^{r}}{r}\,dr.
\label{e_sin_2}
\end{align}
Because $a<b<0$, it follows that $\frac{a}{b} > 1$, and $\frac{\pi
a}{b} > 1$. For $r < 0$, $\frac{\pi a}{b}r < r$, so that
$e^{\frac{\pi a}{b}r}-e^{r} < 0$, and therefore
\begin{equation}
\frac{e^{\frac{\pi a}{b}r}-e^{r}}{r} > 0.
\label{e_sin_2_pos}
\end{equation}
By \ref{e_sin_2_pos}, the integrand in \eqref{e_sin_2} is positive,
so we can bound it with
\begin{align}
\int_{[b,0]}\frac{e^{\frac{\pi a}{b}r}-e^{r}}{r}\,dr
& \leq
\int_{[-\infty,0]}\frac{e^{\frac{\pi a}{b}r}-e^{r}}{r}\,dr
\notag \\
& =
\int_{[-\infty,-1]}\frac{e^{\frac{\pi a}{b}r}-e^{r}}{r}\,dr
+
\int_{[-1,0]}\frac{e^{\frac{\pi a}{b}r}-e^{r}}{r}\,dr
\notag \\
& \leq
\frac{b}{\pi a}e^{-\frac{\pi a}{b}}+e^{-1}
+
\int_{[-1,0]}\frac{e^{\frac{\pi a}{b}r}-e^{r}}{r}\,dr.
\label{e_sin_3}
\end{align}
By Taylor expansion and an error estimate for alternating sums,
\begin{align}
\int_{[-1,0]}\frac{e^{\frac{\pi a}{b}r}-e^{r}}{r}\,dr
& =
\int_{[-1,0]}\sum_{n=1}^{\infty}(\frac{(\frac{\pi a}{b})^n-1}{n!}r^{n-1})\,dr
\notag \\
& =
-\sum_{n=1}^{\infty}(\frac{(\frac{\pi a}{b})^n-1}{n!}\frac{(-1)^{n}}{n})
\notag \\
& \leq
(\frac{\pi a}{b})-1.
\label{e_sin_4}
\end{align}
Combining \eqref{e_sin_3} and \eqref{e_sin_4}, the result follows.
\end{proof}

\begin{lemm}
\label{convol} Let $\varphi(\xi)$ be a function so that the Fourier
transform $\hat{\varphi}(x)$ is a smooth function with support in
$[\frac{1}{2},2]$, and for $\omega \in \mathbb R\setminus \{0\}$,
let $\varphi_{\omega}(\xi)$ be the function such that
$\hat{\varphi}_{\omega}=\varphi(\frac{x}{\omega})$. Then
\begin{equation}
\left |\int_{\mathbb R}\varphi_{\omega}(\xi-z)\frac{e^{itz^3}}{(1+z^2)^{\frac{1}{8}}}\,dz \right |
\lesssim \left\{
\begin{array}{c l}
  (1+t){\omega}^{-\frac{1}{8}} & \textrm{if $\omega >0$,} \\
& \textrm{and }|\xi| \textrm{ intermediate,}\\
(1+t)|\omega|^{-1} & \textrm{else}.
\end{array}
\right.\notag
\end{equation}
\end{lemm}
\begin{proof}
We consider separately the four different cases, ${\omega}<0$,
${\omega}>0$ and $|\xi|$ near, ${\omega}>0$ and $|\xi|$
intermediate, and ${\omega}>0$ and $|\xi|$ far.

\emph{Case ${\omega}<0$:}

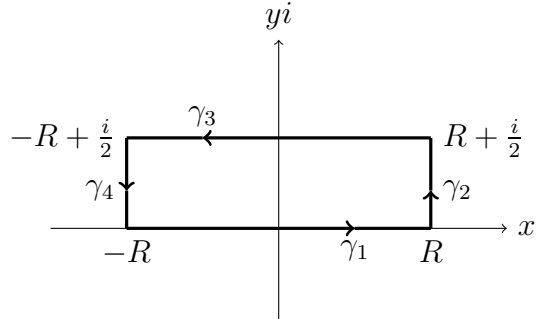
\begin{figure}
    \beginpgfgraphicnamed{fig1}
\begin{center}
\begin{tikzpicture}[domain=0:4]
\draw[->] (-3,0) -- (3,0) node[right] {$x$};

\draw[->] (0,-1.2) -- (0,2.5) node[above] {$yi$};

\draw[very thick,->] (-2,0) --(1,0) node[below] {$\gamma_1$};

\draw[very thick] (1,0) --(2,0) node[below] {$R$};

\draw [very thick, ->] (2,0) -- (2,0.5) node[right] {$\gamma_2$};

\draw [very thick] (2,0.5) -- (2,1.2) node[right] {$R+\frac{i}{2}$};

\draw [very thick, ->] (2,1.2) -- (-1,1.2) node[above] {$\gamma_3$};

\draw [very thick] (-1,1.2) -- (-2,1.2) node[left]
{$-R+\frac{i}{2}$};

\draw [very thick, ->] (-2,1.2) -- (-2,0.5) node[left] {$\gamma_4$};

\draw [very thick] (-2,0.5) -- (-2,0) node[below] {$-R$};
\end{tikzpicture}
\endpgfgraphicnamed
\caption{The contours used for ${\omega} <0$.} \label{Fig1}
\end{center}
\end{figure}

Instead of integrating over $\mathbb R$ in \eqref{project}, we will
compute the integral over the contours $\gamma_1$ through $\gamma_4$
in Figure 1, taking the limit as $R$ approaches infinity. By
Corollary \ref{anal_cont_cor} and the Cauchy integral theorem,
\begin{align}
\int_{\gamma_1}\varphi_{\omega}(\xi-z)\frac{e^{itz^3}}{(1+z^2)^{\frac{1}{8}}}\,dz & =
-\int_{\gamma_2}\varphi_{\omega}(\xi-z)\frac{e^{itz^3}}{(1+z^2)^{\frac{1}{8}}}\,dz
-\int_{\gamma_3}\ldots \notag \\
& \quad -\int_{\gamma_4}\ldots
\notag
\end{align}
We will use estimates on the integrals over $\gamma_2$, $\gamma_3$,
and $\gamma_4$ to estimate \eqref{project}. Along $\gamma_2$,
\begin{align}
& \left |\int_{[0,\frac{1}{2}]}
\varphi_{\omega}(\xi-R-yi)\frac{e^{it(R+yi)^3}}{(1+(R+yi)^2)^{\frac{1}{8}}}i\,dy \right |
\lesssim \notag \\
& \int_{[0,\frac{1}{2}]}
\left | \varphi_{\omega}(\xi-R-yi)\frac{e^{it(R+yi)^3}}{(1+(R+yi)^2)^{\frac{1}{8}}}i \right | \,dy
\lesssim \notag \\
&
\int_{[0,\frac{1}{2}]}
\frac{|e^{2{\omega}y}-e^{\frac{1}{2}{\omega}y}|}{{\omega}^2y|\xi-R-yi|^2}\frac{e^{-t(3R^2-y^2)y}}{(1+R^2)^{\frac{1}{8}}} \,dy. \label{gam_2}
\end{align}
For fixed ${\omega}$, \eqref{gam_2} approaches 0 as $R \rightarrow
\infty$. A similar estimate applies for $\gamma_4$. We can estimate
the integral along $\gamma_3$ using Lemma \ref{anal_cont},
\begin{align}
& \left |\int_{[-R,R]}
\varphi_{\omega}(\xi-x-\frac{i}{2})\frac{e^{it(x+\frac{i}{2})^3}}{(1+(x+\frac{i}{2})^2)^{\frac{1}{8}}}\,dx \right |  \lesssim \notag \\
& \int_{[-R,R]}\left |
\varphi_{\omega}(\xi-x-\frac{i}{2})\frac{e^{it(x+\frac{i}{2})^3}}{(1+(x+\frac{i}{2})^2)^{\frac{1}{8}}}\right | \,dx  \lesssim \notag \\
& \int_{[-R,R]}
\frac{|e^{{\omega}}-e^{\frac{1}{4}{\omega}}|}{{\omega}^2((\xi-x)^2+1)}\frac{e^{-t(\frac{3}{2}x^2-\frac{1}{8})}}{(1+x^2)^{\frac{1}{8}}}\,dx \lesssim \notag \\
& \frac{|e^{{\omega}}-e^{\frac{1}{4}{\omega}}|}{{\omega}^2}\int_{\mathbb R}
\frac{1}{((\xi-x)^2+1)}\frac{1}{(1+x^2)^{\frac{1}{8}}}\,dx \lesssim \notag \\
& \frac{|e^{{\omega}}-e^{\frac{1}{4}{\omega}}|}{{\omega}^2}. \label{gam_3}
\end{align}
From \eqref{gam_2} and \eqref{gam_3} we estimate \eqref{project},
\begin{equation}
\left | \int_{\mathbb R}\varphi_{\omega}(\xi-z)\frac{e^{itz^3}}{(1+z^2)^{\frac{1}{8}}}\,dz \right |
\lesssim \frac{|e^{{\omega}}-e^{\frac{1}{4}{\omega}}|}{{\omega}^2} \lesssim (1+t){|\omega|}^{-1}\notag.
\end{equation}
\emph{End of Case ${\omega}<0$.}

Let $\varepsilon$ be some positive number that will be specified
later. For the remaining three cases, we split up the integral
\eqref{project} in the following manner.
\begin{align}
\int_{\mathbb R}\varphi_{\omega}(\xi-z)\frac{e^{itz^3}}{(1+z^2)^{\frac{1}{8}}}\,dz & =  \int_{\mathbb R \setminus B_{\frac{1}{10}\varepsilon}(\xi)}\varphi_{\omega}(\xi-z)\frac{e^{itz^3}}{(1+z^2)^{\frac{1}{8}}}\,dz
\notag \\
& \quad +
\int_{ B_{\frac{1}{10}\varepsilon}(\xi)}\varphi_{\omega}(\xi-z)\frac{e^{itz^3}}{(1+z^2)^{\frac{1}{8}}}\,dz. \notag
\end{align}
We estimate the integral over $\mathbb R \setminus
B_{\frac{1}{10}\varepsilon}(\xi)$ using the decay of
$\varphi_{\omega}$, from Lemma \ref{anal_cont}.
\begin{align}
\left | \int_{\mathbb R \setminus B_{\frac{1}{10}\varepsilon}(\xi)}\varphi_{\omega}(\xi-z)\frac{e^{itz^3}}{(1+z^2)^{\frac{1}{8}}}\,dz \right | & \leq
\int_{\mathbb R \setminus B_{\frac{1}{10}\varepsilon}(\xi)}\left | \varphi_{\omega}(\xi-z)\frac{e^{itz^3}}{(1+z^2)^{\frac{1}{8}}} \right |\,dz \notag \\
& \leq
\int_{\mathbb R \setminus B_{\frac{1}{10}\varepsilon}(\xi)}\left | \varphi_{\omega}(\xi-z) \right |\,dz
\notag \\
& \leq
\int_{\mathbb R \setminus B_{\frac{1}{10}\varepsilon}(\xi)} \frac{1}{{\omega}(\xi-x)^2} \,dx
\lesssim \frac{1}{{\omega}\varepsilon}. \label{r-ep}
\end{align}
In the next three cases we estimate
\begin{equation}
\int_{ B_{\frac{1}{10}\varepsilon}(\xi)}\varphi_{\omega}(\xi-z)\frac{e^{itz^3}}{(1+z^2)^{\frac{1}{8}}}\,dz.
\label{b+ep}
\end{equation}

\emph{Case $\omega > 0$, near:}
\begin{figure}
    \beginpgfgraphicnamed{fig2}
    \begin{center}
\begin{tikzpicture}[domain=0:4]

\draw[thin] (-2,0) --(-0.5,0);

\draw[thin] (-0.5,0) --(0,0) node[above] {$\xi$};

\draw[thin] (0,0) --(2,0) node[above]
{$\xi+\frac{1}{10}\varepsilon$};

\draw [very thick, ->]+(2,0) arc (0:-45:2) node[right] {$\Gamma_1$};

\draw [very thick]+(1.414,-1.414) arc (-45:-90:2) node[below]
{$\xi-\frac{i}{10}\varepsilon$};

\draw [very thick]+(0,-2) arc (-90:-180:2) node[above]
{$\xi-\frac{1}{10}\varepsilon$};

\end{tikzpicture}
\end{center}
\endpgfgraphicnamed
\caption{The contour used when ${\omega}>0$ and $|\xi| \leq
\frac{1}{10}\sqrt{\frac{{\omega}}{t}}$.} \label{fig2}
\end{figure}
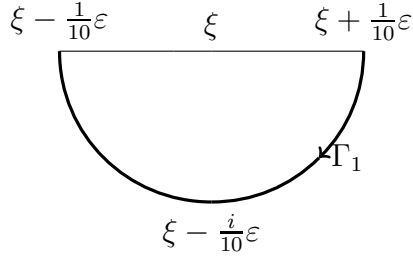

By Corollary \ref{anal_cont_cor} and the Cauchy integral theorem,
we can estimate \eqref{b+ep} by approximating the integral along the
semicircle arc $\Gamma_1$ in Figure 2, as long as we avoid the rays
where the integrand is not analytic. If
$\frac{1}{10}{\omega}^{\frac{1}{2}}t^{-\frac{1}{2}}<1$, then let
$\varepsilon={\omega}^{\frac{1}{2}}t^{-\frac{1}{2}}$. Otherwise, let
$\varepsilon=1$. We illustrate the estimate only for the case
$\varepsilon=1$, as the other case follows by a similar argument.
\begin{align}
& \left | \int_{[2\pi,\pi]}
\varphi_{\omega}(-\frac{\varepsilon}{10} e^{is})
\frac{1}{(1+(\xi+\frac{\varepsilon}{10} e^{is})^2)^{\frac{1}{8}}}
e^{it(\xi+\frac{\varepsilon}{10} e^{is})^3} i\frac{\varepsilon}{10} e^{is}\,ds \right | \lesssim \notag \\
&   \int_{[2\pi, \pi]}
\left | \varphi_{\omega}(-\frac{\varepsilon}{10} e^{is})
\frac{1}{(1+(\xi+\frac{\varepsilon}{10} e^{is})^2)^{\frac{1}{8}}}
e^{it(\xi+\frac{\varepsilon}{10} e^{is})^3} \right | \varepsilon\,ds \lesssim \notag \\
&  \int_{[2\pi,\pi]}t\frac{|e^{\frac{1}{5}{\omega}\varepsilon
\sin{s}}-e^{\frac{1}{20}{\omega}\varepsilon\sin{s}}|}{{\omega}^3\varepsilon\sin{s}}
e^{-\frac{t}{10}(3(\xi+\frac{1}{10}\varepsilon\cos{s})^2-\frac{1}{100}\varepsilon^2\sin^2{s})\varepsilon\sin{s}}
\varepsilon\,ds \label{near_1}
\end{align}
Since $|\xi| \leq \frac{1}{10}\sqrt{\frac{{\omega}}{t}}$ and
$\varepsilon = 1 \leq \sqrt{\frac{{\omega}}{t}}$, it follows that
\begin{align}
|\frac{t}{10}(3(\xi+\frac{1}{10}\varepsilon\cos{s})^2-\frac{1}{100}\varepsilon^2)|
& \leq
\frac{t}{10}(3(|\xi|+\frac{1}{10}\varepsilon)^2+\frac{1}{100}\varepsilon^2)
\notag \\
& \leq \frac{13}{1000}{\omega}.
\notag
\end{align}
Using this and Lemma \ref{e_sin}, we bound \eqref{near_1} with
\begin{align}
 & \int_{[2\pi,\pi]}t\frac{|e^{\frac{1}{5}{\omega}\varepsilon\sin{s}}-e^{\frac{1}{20}{\omega}\varepsilon\sin{s}}|}{{\omega}^3\sin{s}}
e^{-\frac{13}{1000}{\omega}\varepsilon\sin{s}}
\,ds \notag \\
& \quad  \lesssim \frac{t}{{\omega}^3}\int_{[2\pi,\pi]}\frac{|e^{0.187{\omega}\varepsilon\sin{s}}-e^{0.037{\omega}\varepsilon\sin{s}}|}{\sin{s}}
\,ds \lesssim \frac{t}{\omega}. \label{near}
\end{align}
From \eqref{near} and \eqref{r-ep}, we have the estimate
\begin{equation}
\left |\int_{\mathbb R}\varphi_{\omega}(\xi-z)\frac{e^{itz^3}}{(1+z^2)^{\frac{1}{8}}}\,dz \right |
 \lesssim \frac{1}{\omega}+\frac{t}{\omega} \lesssim (1+t){\omega}^{-1}.
\notag
\end{equation}
\emph{End of Case $\omega > 0$, near.}

\emph{Case $\omega > 0$, intermediate:} To estimate \eqref{b+ep}, we use the Young
inequality, and the fact that $\|\varphi_{\omega}\|_1$ is uniformly
bounded in ${\omega}$. Let
$\varepsilon=\frac{1}{10}\sqrt{\frac{{\omega}}{t}}$.
\begin{align}
& \left |
\int_{B_{\frac{1}{10}\varepsilon}(\xi)}
\varphi(\xi-z)\frac{e^{itz^3}}{(1+z^2)^{\frac{1}{8}}}\,dz \right | \notag \\
& \lesssim \left |
\varphi_{\omega} \ast
(\chi_{B_{\frac{1}{10}\varepsilon}(\xi)}(z)
\frac{e^{itz^3}}{(1+z^2)^{\frac{1}{8}}})
\right | \notag \\
& \lesssim \|\varphi_{\omega}\|_1 \left \|\chi_{B_{\frac{1}{10}\varepsilon}(\xi)}(z)
\frac{e^{itz^3}}{(1+z^2)^{\frac{1}{8}}} \right \|_{\infty} \label{inter_1}
\end{align}
Since $\xi$ is intermediate and
$\varepsilon=\frac{1}{10}\sqrt{\frac{{\omega}}{t}}$, any $z \in
B_{\frac{1}{10}\varepsilon}(\xi)$ will obey the estimate $z \approx
\sqrt{\frac{\omega}{t}}$. This estimate on $z$ allows us to bound
the $\|\cdot\|_{\infty}$ term in \eqref{inter_1} by
\begin{equation}
\left \|\chi_{B_{\frac{1}{10}\varepsilon}(\xi)}(z)
\frac{e^{itz^3}}{(1+z^2)^{\frac{1}{8}}} \right \|_{\infty}
\lesssim t^{\frac{1}{8}}|{\omega}|^{-\frac{1}{8}}. \label{inter}
\end{equation}
From \eqref{inter} and \eqref{r-ep}, we have the estimate
\begin{equation}
\left |\int_{\mathbb R}\varphi_{\omega}(\xi-z)\frac{e^{itz^3}}{(1+z^2)^{\frac{1}{8}}}\,dz \right |
 \lesssim \frac{\sqrt{t}}{{\omega}^{\frac{3}{2}}}+t^{\frac{1}{8}}{\omega}^{-\frac{1}{8}} \lesssim (1+t){\omega}^{-\frac{1}{8}}.
\notag
\end{equation}
\emph{End of Case $\omega > 0$, intermediate.}

\emph{Case $\omega > 0$, far:}
\begin{figure}
    \beginpgfgraphicnamed{fig3}
    \begin{center}
\begin{tikzpicture}[domain=0:4]

\draw[thin] (-2,0) --(-0.5,0) ;

\draw[thin] (-0.5,0) --(0,0) node[below] {$\xi$};

\draw[thin] (0,0) --(2,0) node[below]
{$\xi+\frac{1}{10}\varepsilon$};

\draw [very thick, ->]+(2,0) arc (0:45:2) node[above] {$\Gamma_2$};

\draw [very thick]+(1.414,1.414) arc (45:90:2) node[above]
{$\xi+\frac{i}{10}\varepsilon$};

\draw [very thick]+(0,2) arc (90:180:2) node[below]
{$\xi-\frac{1}{10}\varepsilon$};

\end{tikzpicture}
\end{center}
\endpgfgraphicnamed
\caption{The contour used when ${\omega}>0$ and $|\xi| >
10\sqrt{\frac{{\omega}}{t}}$.} \label{fig3}
\end{figure}
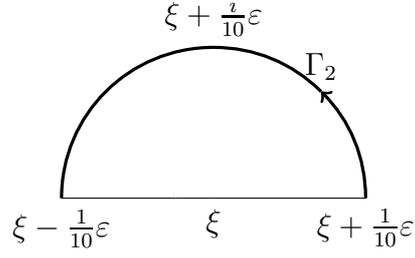
Let $\varepsilon=\sqrt{\frac{{\omega}}{t}}$. We use an argument
similar the near case, integrating along the the semicircle arc
$\Gamma_2$ in Figure 3,
\begin{align}
& \left | \int_{[0,\pi]}
\varphi_{\omega}(-\frac{1}{10}\varepsilon e^{is})
\frac{1}{(1+(\xi+\frac{1}{10}\varepsilon e^{is})^2)^{\frac{1}{8}}}
e^{it(\xi+\frac{1}{10}\varepsilon e^{is})^3} \frac{i}{10}\varepsilon e^{is}\,ds \right | \lesssim \notag \\
&   \int_{[0,\pi]}
\left | \varphi_{\omega}(-\frac{1}{10}\varepsilon e^{is})
\frac{1}{(1+(\xi+\frac{1}{10}\varepsilon e^{is})^2)^{\frac{1}{8}}}
e^{it(\xi+\frac{1}{10}\varepsilon e^{is})^3} \right | \varepsilon\,ds \lesssim \notag \\
&  \int_{[0,\pi]}t\frac{|e^{\frac{1}{5}{\omega}\varepsilon
\sin{s}}-e^{\frac{1}{20}{\omega}\varepsilon \sin{s}}|}{{\omega}^3\varepsilon \sin{s}}
e^{-\frac{t}{10}(3(\xi+\frac{1}{10}\varepsilon \cos{s})^2-\frac{1}{100}\varepsilon^2 \sin^2{s})\varepsilon \sin{s}}
\varepsilon\,ds. \label{far_eq}
\end{align}
Since $\xi > 10\sqrt{\frac{{\omega}}{t}}$,
\begin{align}
& \quad
-29.402{\omega}
\notag \\
& \leq  -\frac{t}{10}(3(10\sqrt{\frac{{\omega}}{t}}-\frac{1}{10}\sqrt{\frac{{\omega}}{t}})^2-\frac{1}{100}\sqrt{\frac{{\omega}}{t}}) \notag \\
& \leq -\frac{t}{10}(3(\xi+\frac{1}{10}\varepsilon \cos{s})^2-\frac{1}{100}\varepsilon^2 \sin^2{s}).
\notag
\end{align}
We use this with Lemma \ref{e_sin} to bound \eqref{far_eq} by
\begin{align}
&  \quad \int_{[0,\pi]}t\frac{|e^{\frac{1}{5}{\omega}\varepsilon \sin{s}}-e^{\frac{1}{20}{\omega}\varepsilon \sin{s}}|}{{\omega}^3\sin{s}}
e^{-29.402{\omega}\varepsilon \sin{s}}
\,ds \notag \\
&
\lesssim \frac{t}{{\omega}^3}\int_{[0,\pi]}\frac{|e^{-29.202{\omega}\varepsilon \sin{s}}-e^{-29.352{\omega}\varepsilon \sin{s}}|}{\sin{s}}
\,ds \lesssim \frac{t}{{\omega}^3}. \label{far}
\end{align}
From \eqref{far} and \eqref{r-ep}, we have the estimate
\begin{equation}
\left |\int_{\mathbb R}\varphi_{\omega}(\xi-z)\frac{e^{itz^3}}{(1+z^2)^{\frac{1}{8}}}\,dz \right |
 \lesssim \frac{\sqrt{t}}{{\omega}^{\frac{3}{2}}}+\frac{t}{{\omega}^3} \lesssim (1+t){\omega}^{-\frac{3}{2}} \lesssim (1+t){\omega}^{-1}.
\notag
\end{equation}
\emph{End of Case $\omega > 0$, far.}
\end{proof}

\begin{lemm}
\label{finite}
\begin{equation}
\left \|D_{\xi}^{\frac{1}{8}} \left (\frac{e^{it\xi^3}}{(1+\xi^2)^{\frac{1}{8}}} \right ) \right \|_{L_{\xi}^{\infty}l_{N}^1}
\lesssim 1+t. \notag
\end{equation}
\end{lemm}
\begin{proof}

The operator $Q^5_N$ (see also Appendix A) is defined by
\begin{equation}
Q^5_Nf \equiv
(|\frac{x}{2^N}|^{\frac{1}{8}}(\eta(\frac{x}{2^N})+\eta(\frac{-x}{2^N}))\hat{f}(x))^{\vee}.
\notag
\end{equation}
Since $Q_N$ is just convolution against the Fourier transform of a
scaled smooth function, by rescaling we obtain
\begin{equation}
\left \|Q_ND_{\xi}^{\frac{1}{8}}\frac{e^{it\xi^3}}{(1+\xi^2)^{\frac{1}{8}}} \right \|_{L_{\xi}^{\infty}l_N^1}
=
\left \|2^{\frac{N}{8}}Q^5_N \left (\frac{e^{it\xi^3}}{(1+\xi^2)^{\frac{1}{8}}} \right ) \right \|_{L_{\xi}^{\infty}l_N^1}.
\notag
\end{equation}

We can estimate the low frequency part using the Young inequality in
the following manner,
\begin{align}
\left \|2^{\frac{N}{8}} Q^5_N \left (\frac{e^{it\xi^3}}{(1+\xi^2)^{\frac{1}{8}}} \right ) \right \|_{L_{\xi}^{\infty}l_{N
\leq 0}^1} & \leq
\sum_{N \leq 0}2^{\frac{N}{8}} \left \|Q^5_N(\frac{e^{it\xi^3}}{(1+\xi^2)^{\frac{1}{8}}}) \right \|_{L_{\xi}^{\infty}} \notag \\
& \lesssim
\sum_{N \leq 0}2^{\frac{N}{8}} \left \|\frac{e^{it\xi^3}}{(1+\xi^2)^{\frac{1}{8}}} \right \|_{L_{\xi}^{\infty}}
\notag \\
& \lesssim \sum_{N \leq 0}2^{\frac{N}{8}} \lesssim 1. \notag
\end{align}
We use Lemma \ref{convol}, noting that if $t$ is fixed, for each
$|\xi|$, there is a unique dyadic $2^N$ so that $\xi$ is
intermediate. We use this to bound the remaining frequencies.
\begin{align}
\left \|2^{\frac{N}{8}}Q^5_N(\frac{e^{it\xi^3}}{(1+\xi^2)^{\frac{1}{8}}}) \right \|_{L_{\xi}^{\infty}l_{N
> 1}^1} & =
\left \|\sum_{N=1}^{\infty}2^{\frac{N}{8}}|Q^5_N(\frac{e^{it\xi^3}}{(1+\xi^2)^{\frac{1}{8}}})| \right \|_{L_{\xi}^{\infty}}
\notag \\
& \lesssim \left \|\sum_{2^N \vert \textrm{ $\xi$ not intermediate}}2^{\frac{N}{8}} \left |Q^5_{N}(\frac{e^{it\xi^3}}{(1+\xi^2)^{\frac{1}{8}}}) \right | \right \|_{L_{\xi}^{\infty}}
\notag \\
& \quad + \left \|\sum_{2^N \vert \textrm{ $\xi$ intermediate}}2^{\frac{N}{8}}|Q^5_{N}(\frac{e^{it\xi^3}}{(1+\xi^2)^{\frac{1}{8}}})| \right \|_{L_{\xi}^{\infty}}
\notag \\
& \lesssim (\sum_{N=1}^{\infty}2^{\frac{N}{8}} 2^{-N} + 1)(1+t). \notag
\end{align}
\end{proof}

\section{Decay Estimates for mKdV Solutions}

With our bound from Lemma \ref{finite}, we will show that our main
result follows. This will come from the fact that for $\alpha \in
(0,1)$,
\begin{equation}
\|D_{x}^{\alpha}(fg)-gD_{x}^{\alpha}f\|_2 \lesssim
\|Q_ND_{x}^{\alpha}g\|_{L_x^{\infty}l_N^1}\|f\|_2.
\label{prod_rule}
\end{equation}
A classical Leibnitz type inequality for fractional derivatives is
the following (see \cite{KPV1}).
\begin{lemm}
\label{stand_prod} Let $0 < \alpha, \alpha_1, \alpha_2 < 1$,
$\alpha=\alpha_1+\alpha_2$, $1 <p,p_1,p_2 < \infty$, and
$\frac{1}{p}=\frac{1}{p_1}+\frac{1}{p_2}$. In addition, the
$\alpha_1=\alpha$, $p=p_2$, and $p_1=\infty$ is allowed. Then the
following holds for functions $f,g$ on $\mathbb R^n$.
\begin{equation}
\|D_x^{\alpha}(fg)-D^{\alpha}(f)g-fD_x^{\alpha}(g)\|_p \lesssim
\|D_x^{\alpha_1}g\|_{p_1}\|D_x^{\alpha_2}f\|_{p_2} \notag
\end{equation}
\end{lemm}
The proof uses the Littlewood-Paley Theorem (see \cite{S}), which
states that for any function $f$, if $1 < p < \infty$, then
\begin{equation}
\label{little_pale}
\|Q_N(f)\|_{L_x^pl_N^2} \lesssim \|f\|_p \lesssim \|Q_N(f)\|_{L_x^pl_N^2}.
\end{equation}

Lemma \ref{stand_prod} is not sufficient for our argument in the
previous section, since we need to put the \emph{derivative} term in
the infinity norm. A product rule like this can be obtained by
following the proof of Lemma \ref{stand_prod} line for line. The
only difference is that since \eqref{little_pale} fails for
$p=\infty$, $\|Q_N(D_x^{\alpha}g)\|_{L_{x}^{\infty}l_N^2}$ is not
equivalent to $\|D_x^{\alpha}g\|_{\infty}$. This idea was inspired
by \cite{KT}, where the authors use
$\|Q_N\cdot\|_{l_N^2L_x^{4}L_T^{\infty}}$ in an estimate where the
$\|\cdot\|_{L_x^{4}L_T^{\infty}}$ norm may fail.
\begin{lemm}
\label{my_prod_rule} Let $0 < \alpha < 1$ and $1 < p < \infty$. For
functions $f$ and $g$,
\begin{equation}
\|D_x^{\alpha}(fg)-gD_x^{\alpha}f-fD_x^{\alpha}g\|_p \lesssim
(\|Q_ND_x^{\alpha}g\|_{L_x^{\infty}l_N^2}
+\|D_x^{\alpha}g\|_{L_x^{\infty}})\|f\|_p. \notag
\end{equation}
In particular,
\begin{equation}
\|D_x^{\alpha}(fg)-gD_x^{\alpha}f\|_2 \lesssim
\|Q_ND_x^{\alpha}g\|_{L_x^{\infty}l_N^1}\|f\|_2. \notag
\end{equation}
\end{lemm}
The proof is in Appendix A.

For a number $1 \leq p \leq \infty$, let $p'$ denote the conjugate
exponent. We recall the following properties of the operator $U(t)$,
\begin{align}
\|\partial_x\int_0^tU(t-t')f(x,t')\,dt'\|_{L_x^2} \lesssim \|f\|_{L_x^1L_T^2}, \label{kato_smooth_dual} \\
%\|\partial_xU(t)f(x)\|_{L_x^{\infty}L_T^2} \lesssim \|f\|_2 \label{kato_smooth} \\
%\|U(t)f(x)\|_{L_x^4L_T^{\infty}} \lesssim \|D_x^{\frac{1}{4}}f\|_2. \label{max_est} \\
\|\int_0^tU(t-t')f(t')\,dt'\|_{L_x^{2}} \lesssim \|f\|_{L_T^{q'}L_x^{p'}}, \label{strichartz} \\
\text{where $p \ge 2$, and $q$ satisfy } \frac{1}{q}=\frac{1}{6}-\frac{1}{3p}. \notag
\end{align}
The proof of \eqref{kato_smooth_dual} can be found in \cite{B}, or
\cite{KPV1}. Inequality \eqref{strichartz} follows from the fact
that $U(t)$ is an $L_x^2$ isometry, along with the dual of the
homogenous Strichartz estimate for $U(t)$ (see \cite{g}, page 1392).

The existence theorem for solutions to \eqref{mkdv} is proved by a
contraction mapping argument, which can also be found in \cite{B}.
\begin{thm}
\label{exist} Let $\|\cdot\|_{Y_T}$ denote the norm such that
\begin{align}
\|f\|_{Y_T} & \equiv \|f\|_{L_x^4L_T^{\infty}}+\|D_x^{\frac{1}{4}}\partial_xf\|_{L_x^{\infty}L_T^{2}} \notag \\
& \quad +\|f\|_{L_T^{\infty}H^{\frac{1}{4}}}+\|\partial_xf\|_{L_x^{20}L_T^{\frac{5}{2}}}+\|D_x^{\frac{1}{4}}f\|_{L_x^{5}L_T^{10}}, \notag \\
 Y_T & \equiv \{f \vert \textrm{ } \|f\|_{Y_T} < \infty\}, \notag
\end{align}
and let $u_0 \in L^2$, and $\Phi$ be the map from $Y_{T}$ to $Y_{T}$
such that
\begin{equation}
\Phi(u) \equiv U(t)u_0-\int_{0}^{t}U(t-t')\partial_x(u^3(t'))\,dt'.
\notag
\end{equation}
Then
\begin{equation}
\|\Phi(u)\|_{Y_T} \lesssim
\|u_0\|_{H^{\frac{1}{4}}}+T^{\frac{1}{2}}\|u\|_{Y_T}^3.
\label{fixed}
\end{equation}
This implies by contraction mapping that there exist
$T=c\|D_x^{\frac{1}{4}}u\|_2^{-4}$ and a unique strong solution
$u(t)$ of the IVP \eqref{mkdv}.
\end{thm}

The proof requires a Leibnitz rule type inequality for
$L_{x}^pL_{T}^q$ norms, which we need as well.
\begin{lemm}
Let $\alpha \in (0,1)$, $\alpha_1,\alpha_2 \in [0,\alpha]$ with
$\alpha=\alpha_1+\alpha_2$. Let $p,q,p_1,p_2,q_2\in (1,\infty)$,
$q_1 \in (1,\infty]$ be such that
\begin{equation}
\frac{1}{p}=\frac{1}{p_1}+\frac{1}{p_2} \textrm{ and } \frac{1}{q}=\frac{1}{q_1}+\frac{1}{q_2}. \notag
\end{equation}
Then
\begin{align}
\|D_x^{\alpha}(fg)-fD_x^{\alpha}g-gD_x^{\alpha}f\|_{L_x^pL_T^q} \lesssim \|D_x^{\alpha_1}f\|_{L_x^{p_1}L_T^{q_1}}\|D_x^{\alpha_2}f\|_{L_x^{p_2}L_T^{q_2}} \notag
\end{align}
Moreover, for $\alpha_1=0$, the value $q_1=\infty$ is allowed.
\end{lemm}

We will need an estimate on the Fourier transform $k(x)$ of
$(1+\xi^2)^{-\frac{1}{8}}$. We expect $k$ to have good decay
properties since it is the inverse Fourier transform of a smooth
function. Since
\begin{equation}
|\hat{\xi}|^{-\frac{1}{4}}=c_0|x|^{-\frac{3}{4}}, \label{basic_ft}
\end{equation}
we expect that $k(x) \approx |x|^{-\frac{3}{4}}$ for small $x$. This
is formalized in the following lemma.
\begin{lemm}
\label{k_est_lem} Let $k(x)$ denote the Fourier transform of the
function $(1+\xi^2)^{-\frac{1}{8}}$. Then for any $n \in \mathbb N$,
\begin{equation}
|k(x)| \lesssim \frac{1}{|x|^{\frac{3}{4}}(1+x^{2n})}. \label{k_est}
\end{equation}
In particular,
\begin{equation}
\int_{\mathbb R}|x|^{\frac{1}{8}}|k(x)| < \infty. \notag
\end{equation}
\end{lemm}
\begin{proof}
For $x>1$, we can repeatedly integrate by parts as follows:
\begin{align}
\int_{\mathbb R}(1+\xi^2)^{-\frac{1}{8}}e^{-ix\xi}\,d\xi
& = \int_{\mathbb R}(1+\xi^2)^{-\frac{1}{8}}\frac{1}{-ix}\frac{d}{\,d\xi}e^{-ix\xi}\,d\xi
\notag \\
& = \frac{1}{ix}\int_{\mathbb R}\frac{1}{4}\xi(1+\xi^2)^{-\frac{9}{8}}e^{-ix\xi}\,d\xi
\notag \\
& = \frac{1}{ix}\int_{\mathbb R}\frac{1}{4}\xi(1+\xi^2)^{-\frac{9}{8}}\frac{1}{-ix}\frac{d}{\,d\xi}e^{-ix\xi}\,d\xi
\notag \\
& = \ldots \notag
\end{align}
This argument gives us the decay in \eqref{k_est}.

When $x < 1$, we split up the integral over the region
$S=[-|x|^{-1},|x|^{-1}]$.
\begin{align}
\int_{\mathbb R}(1+\xi^2)^{-\frac{1}{8}}e^{-ix\xi}\,d\xi
& = \int_{S}(1+\xi^2)^{-\frac{1}{8}}e^{-ix\xi}\,d\xi
\notag \\
& \quad +\int_{\mathbb R \setminus S}(1+\xi^2)^{-\frac{1}{8}}e^{-ix\xi}\,d\xi
\notag \\
& = \mathcal{A}+\mathcal{B}. \notag
\end{align}
Since $(1+\xi^{-2})^{-\frac{1}{8}}$ is bounded,
\begin{align}
|\mathcal{A}| & \lesssim \int_{S}(1+\xi^2)^{-\frac{1}{8}}\,d\xi
\notag \\
& = \int_{S}|\xi|^{-\frac{1}{4}}(1+\xi^{-2})^{-\frac{1}{8}}\,d\xi
\notag \\
& \lesssim \int_{S}|\xi|^{-\frac{1}{4}}\,d\xi
 \lesssim |x|^{-\frac{3}{4}}. \notag
\notag
\end{align}
By integration by parts,
\begin{align}
\mathcal{B} & =  \int_{\mathbb R \setminus S}(1+\xi^2)^{-\frac{1}{8}}\frac{1}{-ix}\frac{d}{\,d\xi}e^{-ix\xi}\,d\xi \notag \\
& = (1+x^{-2})^{-\frac{1}{8}}\frac{e^{i}}{-ix}+(1+x^{-2})^{-\frac{1}{8}}\frac{e^{-i}}{ix} \notag \\
& \quad +\frac{1}{ix}\int_{\mathbb R \setminus S}\frac{1}{4}\xi(1+\xi^2)^{-\frac{9}{8}}e^{-ix\xi}\,d\xi.
\notag
\end{align}
Therefore,
\begin{align}
|\mathcal{B}| & \lesssim  |x|^{-\frac{3}{4}}
+\frac{1}{|x|}\int_{\mathbb R \setminus S}|\xi|(1+\xi^2)^{-\frac{9}{8}}\,d\xi
\notag \\
& =  |x|^{-\frac{3}{4}}
+\frac{1}{|x|}\int_{\mathbb R \setminus S}|\xi|^{-\frac{5}{4}}(1+\xi^{-2})^{-\frac{9}{8}}\,d\xi
\notag \\
& \lesssim  |x|^{-\frac{3}{4}}
+\frac{1}{|x|}\int_{\mathbb R \setminus S}|\xi|^{-\frac{5}{4}}\,d\xi
\lesssim |x|^{-\frac{3}{4}}. \notag
\end{align}
Combining our estimates for $|\mathcal{A}|$ and $|\mathcal{B}|$, the
result follows.
\end{proof}

Before proving Theorem \ref{main}, we prove the corresponding decay
result for solutions to the linear part of \eqref{mkdv}. This is
necessary for the proof of Theorem \ref{main}, and it is also a
simpler case that illustrates the main idea of our proof of Theorem
\ref{main}. We note that it is also possible to prove this result
using an argument like Lemma 2 in \cite{NP}, but this proof does not
generalize to solutions of \eqref{mkdv} as discussed in the
introduction.
\begin{lemm}
\label{main_lin}
For $u_0 \in C_0^{\infty}(\mathbb R)$,
\begin{equation}
\||x|^{s}U(t)u_0(x)\|_2 \lesssim (1+|t|+|t|^{s})\|u_0\|_{H^{2s}}+
\||x|^{s}u_0\|_2. \notag
\end{equation}
\end{lemm}
\begin{proof}
For concreteness, it will suffice to prove the result in the case
$s=\frac{1}{8}$. By the definition of $U(t)$ and the triangle
inequality,
\begin{align} & \||x|^{\frac{1}{8}}
U(t)u_0\|_{2} =
\|D_{\xi}^{\frac{1}{8}}\left ( e^{it\xi^3}\hat{u}_0\right )\|_{L_{\xi}^2} \notag \\
& =\left \| D_{\xi}^{\frac{1}{8}}\left ( \frac{e^{it\xi^3}}{(1+\xi^2)^{\frac{1}{8}}}(1+\xi^2)^{\frac{1}{8}}\hat{u}_0\right )
\right \|_{2} \notag \\
& \lesssim
\left \|D_{\xi}^{\frac{1}{8}}\left ( \frac{e^{it\xi^3}}{(1+\xi^2)^{\frac{1}{8}}}(1+\xi^2)^{\frac{1}{8}}\hat{u}_0\right )
-\frac{e^{it\xi^3}}{(1+\xi^2)^{\frac{1}{8}}}D_{\xi}^{\frac{1}{8}}((1+\xi^2)^{\frac{1}{8}}\hat{u}_0) \right \|_{2} \notag \\
& \quad +\|\frac{e^{it\xi^3}}{(1+\xi^2)^{\frac{1}{8}}}D_{\xi}^{\frac{1}{8}}((1+\xi^2)^{\frac{1}{8}}\hat{u}_0)\|_{2} \notag \\
& \equiv \mathcal{I} + \mathcal{II}. \notag
\end{align}
We can write term $\mathcal{II}$ as
\begin{align}
\mathcal{II}
& = \|(1+\xi^2)^{-\frac{1}{8}}D_{\xi}^{\frac{1}{8}}((1+\xi^2)^{\frac{1}{8}}\hat{u}_0)\|_{2}
\notag \\
& = \|[(1+\xi^2)^{-\frac{1}{8}},D_{\xi}^{\frac{1}{8}}]((1+\xi^2)^{\frac{1}{8}}\hat{u}_0)
+D_{\xi}^{\frac{1}{8}}\hat{u}_0\|_{2}
\notag \\
& \leq \|[(1+\xi^2)^{-\frac{1}{8}},D_{\xi}^{\frac{1}{8}}]((1+\xi^2)^{\frac{1}{8}}\hat{u}_0)\|_2
+\|D_{\xi}^{\frac{1}{8}}\hat{u}_0\|_{2}
\label{term_2}
\end{align}

We need to bound the commutator term in $\mathcal{II}$. For any
function $h$, we use the Plancherel theorem, the Young inequality,
and Lemma \ref{k_est_lem} to obtain
\begin{align}
\|[(1+\xi^2)^{-\frac{1}{8}},D_{\xi}^{\frac{1}{8}}]h\|_{L_{\xi}^2}
& = \|\int_{\mathbb
R}(|x|^{\frac{1}{8}}-|y|^{\frac{1}{8}})k(x-y)\hat{h}(y)\,dy\|_{L_{x}^2}
\notag \\
& \lesssim \|\int_{\mathbb
R}|x-y|^{\frac{1}{8}}|k(x-y)||\hat{h}(y)|\,dy\|_{L_{x}^2}
\notag \\
& \lesssim c_{k}\|\hat{h}\|_{L_{x}^2} = c_{k}\|h\|_{L_{\xi}^2}.
\notag
\end{align}
We apply this to \eqref{term_2},
\begin{equation}
|\mathcal{II}| \lesssim c_k\|u_0\|_{H^{\frac{1}{4}}}+\||x|^{\frac{1}{8}}u_0\|_2. \notag
\end{equation}

For term $\mathcal{I}$, we use Lemma \ref{my_prod_rule} with Lemma
\ref{finite}.
\begin{align}
|\mathcal{I}| & \lesssim \|Q_N\frac{e^{it\xi^3}}{(1+\xi^2)^{\frac{1}{8}}}\|_{L_{x}^{\infty}l_N^1}
\|u_0\|_2. \notag \\
& \lesssim (1+t)\|u_0\|_2. \notag
\end{align}
Combining our estimates for $\mathcal{I}$ and $\mathcal{II}$, the
result follows.
\end{proof}
\begin{proof}[Proof of Theorem \ref{main}]
For concreteness, we prove the result in the most interesting case when
$k=2$, $s=s'=\frac{1}{8}$, and $t>0$. We use a contraction mapping argument
to prove our decay estimate. The resolution space is
\begin{align}
 \|f\|_{Z_T} & \equiv \||x|^{\frac{1}{8}}f\|_{L_T^{\infty}L_x^2} + \|f\|_{Y_T}. \notag \\
 Z_T & \equiv \{f \vert \textrm{ } \|f\|_{Z_T} < \infty\}. \notag
\end{align}

Let $f(t)\equiv \partial_x(u^3(t))$ for convenience, and
consider
\begin{equation}
\Phi(u)(x,t) = U(t)u(x,0)-\int_0^tU(t-t')f(t')dt'. \label{sol_map}
\end{equation}
Multiply \eqref{sol_map} by $|x|^{\frac{1}{8}}$. The
$|x|^{\frac{1}{8}}U(t)u(x,0)$ term is bounded by Lemma
\ref{main_lin} along with a density argument. We concentrate on the
nonlinear term:
\begin{align} & \left \||x|^{\frac{1}{8}}\int_0^t
U(t-t')f(t')\,dt' \right \|_{L_T^{\infty}L_x^2} =
\left \|D_{\xi}^{\frac{1}{8}}\left ( \int_0^te^{i(t-t')\xi^3}f^{\wedge}(t')\,dt' \right ) \right \|_{L_T^{\infty}L_{\xi}^2} \notag \\
& =\left \|D_{\xi}^{\frac{1}{8}}\left (\int_0^t\frac{e^{i(t-t')\xi^3}}{(1+\xi^2)^{\frac{1}{8}}}(1+\xi^2)^{\frac{1}{8}}f^{\wedge}(t') \right )\,dt' \right \|_{{L_T^{\infty}L_{\xi}^2}} \notag \\
& \lesssim
 \|D_{\xi}^{\frac{1}{8}}\left ( \int_0^t\frac{e^{i(t-t')\xi^3}}{(1+\xi^2)^{\frac{1}{8}}}(1+\xi^2)^{\frac{1}{8}}f^{\wedge}(t')\,dt'\right ) \notag \\
& \quad -\int_0^t\frac{e^{i(t-t')\xi^3}}{(1+\xi^2)^{\frac{1}{8}}}D_{\xi}^{\frac{1}{8}}((1+\xi^2)^{\frac{1}{8}}f^{\wedge}(t'))\,dt' \|_{{L_T^{\infty}L_{\xi}^2}} \notag \\
& \quad +\|\int_0^t\frac{e^{i(t-t')\xi^3}}{(1+\xi^2)^{\frac{1}{8}}}D_{\xi}^{\frac{1}{8}}((1+\xi^2)^{\frac{1}{8}}f^{\wedge}(t'))\,dt'\|_{{L_T^{\infty}L_{\xi}^2}} \notag \\
& \equiv I + II. \label{def-I}
\end{align}
We bound term $II$ in a similar fashion to term $\mathcal{II}$ in
Lemma \ref{main_lin}:
\begin{align}
II & \lesssim \|\int_0^t {e^{i(t-t')\xi^3}}[\frac{1}{(1+\xi^2)^{\frac{1}{8}}},D_{\xi}^{\frac{1}{8}}]((1+\xi^2)^{\frac{1}{8}}f^{\wedge}(t'))\,dt'\|_{{L_T^{\infty}L_{\xi}^2}} \notag \\
& \quad + \|\int_0^t {e^{i(t-t')\xi^3}}D_{\xi}^{\frac{1}{8}}(f^{\wedge}(t'))\,dt'\|_{L_T^{\infty}L^2_{\xi}} \notag \\
& \lesssim \|[\frac{1}{(1+\xi^2)^{\frac{1}{8}}},D_{\xi}^{\frac{1}{8}}]((1+\xi^2)^{\frac{1}{8}}f^{\wedge}(t'))\|_{L_T^1L^2_{\xi}} \notag \\
& \quad + \|\int_0^t {e^{i(t-t')\xi^3}}D_{\xi}^{\frac{1}{8}}(f^{\wedge}(t'))\,dt'\|_{{L_T^{\infty}L_{\xi}^2}} \notag \\
& \lesssim \|((1+\xi^2)^{\frac{1}{8}}f^{\wedge}(t'))\|_{L_T^1L^2_{\xi}} + \|\int_0^t {e^{i(t-t')\xi^3}}D_{\xi}^{\frac{1}{8}}(f^{\wedge}(t'))\,dt'\|_{{L_T^{\infty}L_{\xi}^2}} \notag \\
& \lesssim \|((1+D_x^2)^{\frac{1}{8}}f(t'))\|_{L_T^1L^2_{x}} + \|\int_0^{t}U(t-t')|x|^{\frac{1}{8}}f(t')\,dt'\|_{{L_T^{\infty}L_x^2}} \notag \\
& \equiv II.1+II.2. \notag
\end{align}
Specializing to the case of the mKdV,
$f(t')=\partial_{x}(u^3(t'))$, we bound $II.1$ using
Theorem \ref{exist}:
\begin{align}
II.1 & \lesssim\|\partial_x(u^3)\|_{L_T^1L_x^2} +\|D_x^{\frac{1}{4}}\partial_x(u^3)\|_{L_T^1L_x^2} \notag \\
& \lesssim
T^{\frac{1}{2}}\|\partial_x(u^3)\|_{L_T^2L_x^2} +T^{\frac{1}{2}}\|D_x^{\frac{1}{4}}\partial_x(u^3)\|_{L_T^2L_x^2}
\notag \\
&\lesssim
T^{\frac{1}{2}}\|u\|_{L_x^{4}L_T^{\infty}}^2\|\partial_xu\|_{L_x^{\infty}L_T^{2}}+T^{\frac{1}{2}}\|u^2\|_{L_x^2L_T^{\infty}}\|D_x^{\frac{1}{4}}\partial_xu\|_{L_x^{\infty}L_T^2}
\notag \\
& \quad
+T^{\frac{1}{2}}\|D_x^{\frac{1}{4}}(u^2)\|_{L_x^{\frac{20}{9}}L_T^{10}}\|\partial_xu\|_{L_x^{20}L_T^{\frac{5}{2}}}
\notag \\
& \lesssim T^{\frac{1}{2}}\| u\|_{L_x^{4}L_T^{\infty}}^2\|\partial_xu\|_{L_x^{\infty}L_T^{2}}+
T^{\frac{1}{2}}\|u\|_{L_x^4L_T^{\infty}}^2\|D_x^{\frac{1}{4}}\partial_xu\|_{L_x^{\infty}L_T^2}
\notag \\
& \quad
+T^{\frac{1}{2}}\|u\|_{L_x^{4}L_T^{\infty}}\|D_x^{\frac{1}{4}}u\|_{L_x^{5}L_T^{10}}\|\partial_xu\|_{L_x^{20}L_T^{\frac{5}{2}}}
\notag \\
& \lesssim T^{\frac{1}{2}}\|u\|_{Z_T}^3. \notag
\end{align}
Let $\phi(x) \in C_0^{\infty}(\mathbb R)$ have the property that
$\phi(x)=1$ for $x \in (-1,1)$. We handle $II.2$ with the following
argument:
\begin{align}
&
\|\int_0^{t}U(t-t')|x|^{\frac{1}{8}}f(t')\,dt'\|_{L_T^{\infty}L_x^2}
\lesssim
\|\int_0^{t}U(t-t')|x|^{\frac{1}{8}}\phi(x)f(t')\,dt'\|_{L_T^{\infty}L_x^2} \notag \\
& \quad +
\|\int_0^{t}U(t-t')|x|^{\frac{1}{8}}(1-\phi(x))f(t')\,dt'\|_{L_T^{\infty}L_x^2} \notag \\
& \lesssim \|\int_0^{t}U(t-t')[|x|^{\frac{1}{8}}(1-\phi(x)),\partial_x]u^3(t')\,dt'\|_{L_T^{\infty}L_x^2} \notag \\
& \quad +\|\int_0^{t}U(t-t')\partial_x((|x|^{\frac{1}{8}}(1-\phi(x)))u^3(t'))\,dt'\|_{L_T^{\infty}L_x^2}\notag \\
& \quad +\|\int_0^{t}U(t-t')|x|^{\frac{1}{8}}\phi(x)\partial_x(u^3(t'))\,dt'\|_{L_T^{\infty}L_x^2}\notag \\
& \equiv II.2.a+II.2.b+II.2.c. \notag
\end{align}
For $II.2.a$, we use \eqref{strichartz}, and that for any function
$h$, and $p \ge 1$,
\begin{equation}
\|[|x|^{\frac{1}{8}}(1-\phi(x)),\partial_x]h\|_p \lesssim
\|\frac{\partial}{\partial
x}(|x|^{\frac{1}{8}}(1-\phi(x)))\|_{\infty}\|h\|_p, \notag
\end{equation}
along with the Sobolev inequality to obtain the bound
\begin{align}
II.2.a & \lesssim \|u^3\|_{L_T^{\frac{12}{11}}L_x^{\frac{4}{3}}}=\|u\|_{L_T^{\frac{36}{11}}L_x^{4}}^3 \notag \\
& \lesssim T^{\frac{11}{12}}\|u\|_{L_{T}^{\infty}H^{\frac{1}{4}}}^3 \lesssim T^{\frac{11}{12}}\|u\|_{Z_T}^3. \notag
\end{align}
We use \eqref{kato_smooth_dual} to estimate $II.2.b$:
\begin{align}
II.2.b & \lesssim \||x|^{\frac{1}{8}}(1-\phi(x))u^3\|_{L_x^1L_T^2} \notag \\
& \lesssim \|u^2\|_{L_x^2L_T^{\infty}}\||x|^{\frac{1}{8}}(1-\phi(x))u\|_{L_x^2L_T^{2}} \notag \\
& \lesssim
\|u\|_{L_x^4L_T^{\infty}}^2(\||x|^{\frac{1}{8}}u\|_{L_T^{2}L_x^2}
+\|u\|_{L_T^{2}L_x^2}) \notag \\
& \lesssim
T^{\frac{1}{2}}\|u\|_{L_x^4L_T^{\infty}}^2(\||x|^{\frac{1}{8}}u\|_{L_T^{\infty}L_x^2}
+\|u\|_{L_T^{\infty}L_x^2}) \notag \\
&  \lesssim
T^{\frac{1}{2}}\|u\|_{Z_T}^3. \notag
\end{align}
We use Theorem \ref{exist} and the fact that $\phi$ has compact
support to control $II.2.c$:
\begin{align}
II.2.c & \lesssim \||x|^{\frac{1}{8}}\phi(x)\partial_x(u^3)\|_{L_T^1L_x^2} \notag \\
& \lesssim \|\partial_x(u^3)\|_{L_T^1L_x^2} \notag \\
& \lesssim T^{\frac{1}{2}}\|u\|_{L_x^4L_T^{\infty}}^2\|\partial_xu\|_{L_x^{\infty}L_T^{2}} \notag \\
& \lesssim T^{\frac{1}{2}}\|u\|_{Z_T}^3. \notag
\end{align}
Term $I$ from \eqref{def-I} can be controlled using Theorem
\ref{exist}, and the same argument as the bound for $II.1$:
\begin{align}
I & \lesssim \|Q_ND_{\xi}^{\frac{1}{8}}\frac{e^{it\xi^3}}{(1+\xi^2)^{\frac{1}{8}}}\|_{L_T^{2}L_{\xi}^{\infty}l_N^1}
\|(1+D_x^2)^{\frac{1}{8}}(u^2\partial_x u)\|_{L_T^2L_{x}^2} \notag \\
& \lesssim (1+T^{\frac{3}{2}})
(\|u^2\partial_x u\|_{L_T^2L_{x}^2}+\|D_x^{\frac{1}{4}}(u^2\partial_x u)\|_{L_T^2L_{x}^2}) \notag \\
& \lesssim
T^{\frac{1}{2}}(1+T)
(\|u\|_{L_{x}^{4}L_{T}^{\infty}}^2\|\partial_xu\|_{L_{x}^{\infty}L_{T}^{2}}
+ \|u\|_{L_x^4L_T^{\infty}}^2\|D_x^{\frac{1}{4}}\partial_xu\|_{L_x^{\infty}L_T^2} \notag \\
& \quad +
\|u\|_{L_x^{4}L_T^{\infty}}\|D_x^{\frac{1}{4}}u\|_{L_x^{5}L_T^{10}}\|\partial_xu\|_{L_x^{20}L_T^{\frac{5}{2}}})
\notag \\
& \lesssim T^{\frac{1}{2}}(1+T)\|u\|_{Z_T}^3
\lesssim T^{\frac{1}{2}}(1+T)\|u\|_{Z_T}^3. \notag
\end{align}
Putting these estimates together,
\begin{align}
\||x|^{\frac{1}{8}}u\|_{L_T^{\infty}L_x^2} & \lesssim
\||x|^{\frac{1}{8}}U(t)u_0\|_{L_T^{\infty}L_x^2}+I+II.1+II.2.a+II.2.b+II.2.c \notag \\
& \lesssim \||x|^{\frac{1}{8}}u_0\|_{L_x^2}+(1+T)\|u_0\|_{H^{\frac{1}{4}}}
+T^{\frac{1}{2}}(1+T^{\frac{5}{12}}+T)\|u\|_{Z_T}^3. \label{contract}
\end{align}
In order to get a contraction, we need to bound $\|u\|_{Y_T}$ in
terms of $\|u\|_{Z_T}$. This follows from estimate \eqref{fixed} in
Theorem \ref{exist}. By combining this with \eqref{contract}, we
obtain a contraction by taking $T$ small enough,
\begin{equation}
\|u\|_{Z_T} \lesssim
\||x|^{\frac{1}{8}}u_0\|_{L_x^2}+(1+T)\|u_0\|_{H^{\frac{1}{4}}}
+T^{\frac{1}{2}}(1+T^{\frac{5}{12}}+T)\|u\|_{Z_T}^3. \notag
\end{equation}
In order to show that $\||x|^{\frac{1}{8}}u(t)\|_{L_x^2}$ is finite,
for $t \in [0,T)$, apply $\||x|^{\frac{1}{8}}\cdot\|_2$ to
\eqref{sol_map} instead of
$\||x|^{\frac{1}{8}}\cdot\|_{L_T^{\infty}L_x^2}$, keeping in mind
that $\||x|^{\frac{1}{8}}u\|_{L_T^{\infty}L_x^2}$ is finite.
\end{proof}

\section{Appendix A}
For our proof of Lemma \ref{prod_rule}, we closely follow the proof
of Theorem A.8 in \cite{KPV1}. This requires more notation. Let
$\alpha_1=0$, $\alpha_2=\alpha \in [0,1]$. For a function $f$, let
\begin{equation}
P_Nf \equiv \sum_{j \leq N-3}Q_jf. \notag
\end{equation}
Define $p(x)$ to be the function so that
\begin{equation}
(P_Nf)^{\wedge} = p(2^{-N}x)\hat{f}. \notag
\end{equation}
Let $\tilde{p} \in C_{0}^{\infty}(\mathbb R)$, with $\tilde{p}(x)
=1$  for $x \in [-100,100]$, and let
\begin{equation}
(\tilde{P}_Nf)^{\wedge}(x)=\tilde{p}(2^{-N}x)\hat{f}.
\notag
\end{equation}
Let $\tilde{\eta}\in C_0^{\infty}(\mathbb R)$ with
$\tilde{\eta}(x)=1$ for $x \in [\frac{1}{4},4]$, and
supp$\tilde{\eta} \in [\frac{1}{8},8]$. Then define
$(\tilde{Q}_kf)^{\wedge}(x)=\tilde{\eta}(2^{-k}x)\hat{f}$. Let
\begin{align}
\Psi^{i}(x)=|x|^{\alpha_j}p(x)
& \textrm{, }
\eta^{j}(x)=\frac{\eta(x)}{|x|^{\alpha_j}},
\notag \\
(\Psi_j^k)^{\wedge}(x)=\Psi_j(2^{-k})\hat{f}(x)
& \textrm{, and }
(Q_k^jf)^{\wedge}(x)=\eta^j(2^{-k}x)\hat{f}(x).
\notag
\end{align}
Similarly, with $\eta^3(x)=|x|^{\alpha}\tilde{p}(x)$, $\eta^4(x)=|x|^{\alpha_1}\eta(x)$, and $\eta^5(x)=|x|^{\alpha_2}\eta(x)$ we define $Q_k^3,Q_k^4,Q_k^5$.
Let
\begin{align}
\eta^{\nu,j}(x) & =\exp(i\nu x)\eta^j(x), \notag \\
\eta^{\mu,j}(x) & =\exp(i\mu x)x|x|^{-\alpha_j}p(x), \notag
\end{align}
with $j=1,2$ and $Q_k^{\nu,j}, Q_k^{\mu,j}$ the corresponding operators.

The following is Proposition A.2 from \cite{KPV1}.
\begin{lemm}
\label{prop_a2}
\begin{align}
D_x^{\alpha}(fg)-&fD_x^{\alpha}g-gD_x^{\alpha}f
\notag \\
&=
\sum_{|j|<2}2^{j\alpha_2}\sum_{k}Q_k^3(Q_k^1(D^{\alpha_1}f)Q_{k-j}^2(D^{\alpha_2}g))
\notag \\
& \quad + \sum_{k}\tilde{Q}_k(\Psi_k^1(D^{\alpha_2}g)Q_k^1(D^{\alpha_1}f))
\notag \\
& \quad + \sum_{k}\tilde{Q}_k(Q_k^2(D^{\alpha_2}g)\Psi_k^2(D^{\alpha_1}f))
\notag \\
& \quad + \sum_{|j| \leq 2}2^{j\alpha_2}\sum_{k}Q_k^1(D^{\alpha_1}f)Q_{k-j}^4(D^{\alpha_2}g)
\notag \\
& \quad + \sum_{|j| \leq 2}2^{j\alpha_2}\sum_{k}Q_{k-j}^2(D^{\alpha_2}g)Q_{k}^5(D^{\alpha_1}f)
\notag \\
& \quad +\int_{\mathbb R}\int_{\mathbb R}
\left [
\sum_{k} \tilde{Q}_k(Q_k^{\nu,1}(D^{\alpha_1}f)Q_k^{\mu,2}(D^{\alpha_2}g))
\right ]r_1(\mu,\nu)\,d\nu\,d\mu
\notag \\
& \quad +\int_{\mathbb R}\int_{\mathbb R}
\left [
\sum_{k} \tilde{Q}_k(Q_k^{\nu,2}(D^{\alpha_2}g)Q_k^{\mu,1}(D^{\alpha_1}f))
\right ]r_2(\mu,\nu)\,d\nu\,d\mu,
\notag
\end{align}
where $r_1,r_2 \in \mathcal{S}(\mathbb R^2)$.
\end{lemm}
\begin{proof}[Proof of Lemma \ref{my_prod_rule}]
From Lemma \ref{prop_a2}, we need to bound four types of terms:
\begin{enumerate}
\item $\sum_{-\infty}^{\infty} Q_k(Q_k(f)Q_k(D_x^{\alpha}g))$

\item $\sum_{-\infty}^{\infty} Q_k(\Psi_k(f)Q_k(D_x^{\alpha}g))$

\item $\sum_{-\infty}^{\infty} Q_k(Q_k(f)\Psi_k(D_x^{\alpha}g))$

\item $\sum_{-\infty}^{\infty} Q_k(f)Q_k(D_x^{\alpha}g)$
\end{enumerate}

Let $\mathcal{M}h$ denote the Hardy Maximal operator applied to the
function $h$. We control the first term using duality,
\begin{align}
& |\int_{\mathbb R} \sum_{-\infty}^{\infty} Q_k(Q_k(f)Q_k(D_x^{\alpha}g))h\,dx| = |\int_{\mathbb R} \sum_{-\infty}^{\infty} Q_k(f)Q_k(D_x^{\alpha}g)Q_k(h)\,dx|
\notag \\
& \lesssim \int_{\mathbb R} \sqrt{\sum_{-\infty}^{\infty} |Q_k(f)|^2|Q_k(D_x^{\alpha}g)|^2}\sqrt{\sum_{-\infty}^{\infty}|Q_n(h)|^2}\,dx
\notag \\
& \lesssim \|Q_k(f)Q_k(D_x^{\alpha}g)\|_{L_x^pl_k^2}\|Q_n(h)\|_{L_x^{p'}l_n^2}
\notag \\
& \lesssim \|\mathcal{M}(f)\|_{L_x^p}\|Q_k(D_x^{\alpha}g)\|_{L_x^{\infty}l_k^2}\|Q_n(h)\|_{L_x^{p'}l_n^2}
\notag \\
& \lesssim \|f\|_{p'}\|Q_k(D_x^{\alpha}g)\|_{L_x^{\infty}l_k^2}\|Q_n(h)\|_{L_x^{p'}l_n^2}
\notag \\
& \lesssim
\|f\|_p\|Q_k(D_x^{\alpha}g)\|_{L_x^{\infty}l_k^2}\|h\|_{L_x^{p'}}.
\notag
\end{align}

The second item is treated as the first, with $\Psi_k(f)$ replacing
$Q_k(f)$. A similar argument is used on the third term, with
$\Psi_k(D_x^{\alpha}g)$ replacing $\Psi_k(f)$, and the fact that
\begin{equation}
 \|\mathcal{M}(D_x^{\alpha}g)\|_{L_x^{\infty}} \lesssim \|g\|_{L_x^{\infty}}, \notag
\end{equation}
because $\mathcal{M}$ is a bounded operator from $L^{\infty}$ to
$L^{\infty}$.

The Last term is treated with Cauchy-Schwartz,
\begin{align}
\|\sum_{-\infty}^{\infty}Q_k(f)Q_k(D_x^{\alpha}g)\|_{p} & \lesssim \|\|Q_n(f)\|_{l_n^2}\|Q_k(D_x^{\alpha}g)\|_{l_k^2}\|_p \notag \\
& \lesssim
\|Q_n(f)\|_{L_x^{p}l_k^2}\|Q_k(D_x^{\alpha}g)\|_{L_x^{\infty}l_k^2}.
\notag
\end{align}
This proves the the first part of the lemma.

The second part follows from
\begin{equation}
\|D_x^{\alpha}(fg)-fD_x^{\alpha}g-gD_x^{\alpha}f\|_p \ge
\|D_x^{\alpha}(fg)-gD_x^{\alpha}f\|_p-\|fD_x^{\alpha}g\|_p,
\notag
\end{equation}
the observation that
\begin{equation}
|D^{\alpha}g| \leq \sum_{N}|Q_N(D^{\alpha}g)|,
\notag
\end{equation}
and for arbitrary functions $\varphi_N$,
\begin{equation}
\|\varphi_N\|_{l_N^2} \leq \|\varphi_N\|_{l_N^{\infty}}^{\frac{1}{2}}\|\varphi_N\|_{l_N^1}^{\frac{1}{2}} \leq \|\varphi_N\|_{l_N^1}. \notag
\end{equation}
\end{proof}

\textbf{Acknowledgments:} I would like to thank Gustavo Ponce for many fruitful discussions, and Luiz Farah for reading earlier drafts.


\begin{thebibliography}{10}

\bibitem{MR0079154}
I.~Bihari.
\newblock A generalization of a lemma of {B}ellman and its application to
  uniqueness problems of differential equations.
\newblock {\em Acta Math. Acad. Sci. Hungar.}, 7:81--94, 1956.

\bibitem{MR2018661}
Michael Christ, James Colliander, and Terrence Tao.
\newblock Asymptotics, frequency modulation, and low regularity ill-posedness
  for canonical defocusing equations.
\newblock {\em Amer. J. Math.}, 125(6):1235--1293, 2003.

\bibitem{CKSTT}
J.~Colliander, M.~Keel, G.~Staffilani, H.~Takaoka, and T.~Tao.
\newblock Sharp global well-posedness for {K}d{V} and modified {K}d{V} on
  {$\Bbb R$} and {$\Bbb T$}.
\newblock {\em J. Amer. Math. Soc.}, 16(3):705--749 (electronic), 2003.

\bibitem{PonceFons}
Germ\'{a}n Fonseca and Gustavo Ponce.
\newblock The {I}{V}{P} for the {B}enjamin {O}no equation in weighted sobolev
  spaces.
\newblock {\em arXiv:1004.5592v2}, 2010.

\bibitem{PhysRevLett.19.1095}
Clifford~S. Gardner, John~M. Greene, Martin~D. Kruskal, and Robert~M. Miura.
\newblock Method for solving the korteweg-devries equation.
\newblock {\em Phys. Rev. Lett.}, 19(19):1095--1097, Nov 1967.

\bibitem{g}
J.~Ginibre and Y.~Tsutsumi.
\newblock Uniqueness of solutions for the generalized {K}orteweg-de {V}ries
  equation.
\newblock {\em SIAM J. Math. Anal.}, 20(6):1388--1425, 1989.

\bibitem{MR2531556}
Zihua Guo.
\newblock Global well-posedness of {K}orteweg-de {V}ries equation in
  {$H^{-3/4}(\Bbb R)$}.
\newblock {\em J. Math. Pures Appl. (9)}, 91(6):583--597, 2009.

\bibitem{MR847012}
Nakao Hayashi, Kuniaki Nakamitsu, and Masayoshi Tsutsumi.
\newblock On solutions of the initial value problem for the nonlinear
  {S}chr\"odinger equations in one space dimension.
\newblock {\em Math. Z.}, 192(4):637--650, 1986.

\bibitem{MR880978}
Nakao Hayashi, Kuniaki Nakamitsu, and Masayoshi Tsutsumi.
\newblock On solutions of the initial value problem for the nonlinear
  {S}chr\"odinger equations.
\newblock {\em J. Funct. Anal.}, 71(2):218--245, 1987.

\bibitem{MR987792}
Nakao Hayashi, Kuniaki Nakamitsu, and Masayoshi Tsutsumi.
\newblock Nonlinear {S}chr\"odinger equations in weighted {S}obolev spaces.
\newblock {\em Funkcial. Ekvac.}, 31(3):363--381, 1988.

\bibitem{Ka}
Tosio Kato.
\newblock On the {C}auchy problem for the (generalized) {K}orteweg-de {V}ries
  equation.
\newblock In {\em Studies in applied mathematics}, volume~8 of {\em Adv. Math.
  Suppl. Stud.}, pages 93--128. Academic Press, New York, 1983.

\bibitem{MR1230283}
Carlos Kenig, Gustavo Ponce, and Luis Vega.
\newblock The {C}auchy problem for the {K}orteweg-de {V}ries equation in
  {S}obolev spaces of negative indices.
\newblock {\em Duke Math. J.}, 71(1):1--21, 1993.

\bibitem{KPV1}
Carlos Kenig, Gustavo Ponce, and Luis Vega.
\newblock Well-posedness and scattering results for the generalized korteweg-de
  vries equation via the contraction principle.
\newblock {\em Comm. Pure Appl. Math.}, 46(4):527--620, 1993.

\bibitem{KT}
Carlos Kenig and Hideo Takaoka.
\newblock Global wellposedness of the modified {B}enjamin-{O}no equation with
  initial data in {$H^{1/2}$}.
\newblock {\em Int. Math. Res. Not.}, pages Art. ID 95702, 44, 2006.

\bibitem{Kish}
Nobu Kishimoto.
\newblock Well-posedness of the {C}auchy problem for the {K}orteweg-de {V}ries
  equation at the critical regularity.
\newblock {\em Differential Integral Equations}, 22(5-6):447--464, 2009.

\bibitem{2009ChPhB..18.4074L}
{Z.-L.} {Li}.
\newblock {Application of higher-order KdV-mKdV model with higher-degree
  nonlinear terms to gravity waves in atmosphere}.
\newblock {\em Chinese Physics B}, 18:4074--4082, October 2009.

\bibitem{B}
F.~Linares and G.~Ponce.
\newblock {\em Introduction to Nonlinear Dispersive Equations}.
\newblock Universitext. Springer, New York, 2009.

\bibitem{NP}
J.~Nahas and G.~Ponce.
\newblock On the persistent properties of solutions to semi-linear
  {S}chr\"odinger equation.
\newblock {\em Comm. Partial Differential Equations}, 34(10-12):1208--1227,
  2009.

\bibitem{PRUD}
V.~V. Prudskikh.
\newblock Ion-acoustic solitons in bi-ion dusty plasma.
\newblock {\em Plasma Physics Reports}, 34(11):955--962, Nov. 2008.

\bibitem{1994JNS.....4..355R}
E.~A. {Ralph} and L.~{Pratt}.
\newblock {Predicting eddy detachment for an equivalent barotropic thin jet}.
\newblock {\em Journal of NonLinear Science}, 4:355--374, December 1994.

\bibitem{S}
Elias~M. Stein.
\newblock {\em Singular integrals and differentiability properties of
  functions}.
\newblock Princeton Mathematical Series, No. 30. Princeton University Press,
  Princeton, N.J., 1970.

\bibitem{MR2286393}
Terence Tao.
\newblock Scattering for the quartic generalised {K}orteweg-de {V}ries
  equation.
\newblock {\em J. Differential Equations}, 232(2):623--651, 2007.

\end{thebibliography}
\end{document}